\documentclass[a4paper,reqno]{amsart}
      \usepackage{etex}
      \usepackage[OT2, T1]{fontenc}
      \usepackage{url}
      \usepackage{amsmath}
      \usepackage{lmodern}
      \usepackage{array}
      \usepackage{graphicx}
      \usepackage{amsfonts}
      \usepackage{amssymb}
      \usepackage{extpfeil}
      \usepackage{mathtools}
\usepackage{soul}
\usepackage{xargs}
      \usepackage{xcolor}
      \definecolor{imperialred}{RGB}{237, 41, 57}
      \definecolor{royalblue}{RGB}{64, 106, 212}
      \definecolor{link}{RGB}{11,0,128}
      \definecolor{olivegreen}{RGB}{128, 128, 0}
      \usepackage{amsthm}
      \usepackage[alphabetic,lite]{amsrefs}  % for bibliography
      \usepackage{esint}
      \usepackage{epic}
      \usepackage{calligra}
      \usepackage{xkeyval}

      \usepackage{exscale, relsize}
      \usepackage{stackengine}
      \usepackage{scalerel}
      \usepackage{setspace}
      \usepackage{paralist}
      \usepackage{colonequals}		% for nice :=
      \usepackage[left]{lineno}
      \usepackage{amscd}
      \usepackage[all,cmtip]{xy}	% for commutative diagrams
      \usepackage{sseq}			% for spectral sequences
      \usepackage{verbatim}
      \usepackage{parskip}
      \usepackage{microtype}
      \usepackage{ragged2e}
      \usepackage[in]{fullpage}
      \usepackage{graphicx}
      \usepackage{mathrsfs} 			% for \mathscr (script letters)
      \usepackage{bm} 				% for bold Greek letters
      \usepackage{enumitem}
      \usepackage{moreenum}			% for enumeration via Greek letters
      \usepackage{blindtext}
      \usepackage{tikz-cd}
      \usepackage{tikz}
      \usepackage{subfiles}
      \usetikzlibrary{shapes.geometric}
      \usetikzlibrary{automata, bending,matrix,arrows,decorations.pathmorphing,backgrounds,positioning,fit,petri,arrows.meta}
      \usepackage{textcomp}
      \usepackage{indentfirst}
      \usetikzlibrary{arrows}
      \usepackage{enumitem}
      \tikzset{commutative diagrams/.cd,arrow style=tikz,diagrams={>=latex'}}
      \usepackage{stmaryrd}	% this is for \llbracket and \rrbracket
      \usepackage{subfiles}
      \usepackage{ragged2e}
      \usetikzlibrary{arrows}
      \usepackage{stackengine}
\usepackage[disable, colorinlistoftodos,prependcaption,textsize=tiny]{todonotes}
\usepackage[colorlinks=true, citecolor=royalblue,linkcolor=olivegreen,urlcolor=link]{hyperref}
       \usepackage{cleveref}
\newcommandx{\unsure}[2][1=]{\todo[linecolor=red,backgroundcolor=red!25,bordercolor=red,#1]{#2}}
\newcommandx{\change}[2][1=]{\todo[linecolor=blue,backgroundcolor=blue!25,bordercolor=blue,#1]{#2}}
\newcommandx{\info}[2][1=]{\todo[linecolor=OliveGreen,backgroundcolor=OliveGreen!25,bordercolor=OliveGreen,#1]{#2}}
\newcommandx{\improvement}[2][1=]{\todo[linecolor=Plum,backgroundcolor=Plum!25,bordercolor=Plum,#1]{#2}}
\newcommandx{\thiswillnotshow}[2][1=]{\todo[disable,#1]{#2}}

\usepackage{marginnote}
\newcommand{\mytodo}[2][]{{%
 \let\marginpar\marginnote
 \reversemarginpar
 \renewcommand{\baselinestretch}{0.8}%
 \todo[#1]{#2}}}

         \newcommand{\GG}{\Gamma}
         
         \newcommand{\GGL}{\Lambda}

         \newcommand{\bA}{\mathbb{A}}
         \newcommand{\bC}{\mathbb{C}}

         \newcommand{\bG}{\mathbb{G}}

         \newcommand{\bZ}{\mathbb{Z}}

         \newcommand{\cA}{\mathcal{A}}
         \newcommand{\cB}{\mathcal{B}}

         \newcommand{\cI}{\mathcal{I}}

         \newcommand{\fm}{\mathfrak{m}}
         
         \newcommand{\fp}{\mathfrak{p}}

         \newcommand{\sE}{\mathscr{E}}
         \newcommand{\sF}{\mathscr{F}}
         
         \newcommand{\sH}{\mathscr{H}}

         \newcommand{\sL}{\mathscr{L}}
         \newcommand{\sM}{\mathscr{M}}
         \newcommand{\sN}{\mathscr{N}}
         \newcommand{\sO}{\mathscr{O}}

         \newcommand{\ra}{\rightarrow}

         \newcommand{\Lra}{\Leftrightarrow}

         \newcommand{\lra}{\longrightarrow}
         
         \newcommand{\hra}{\hookrightarrow}

         \newcommand{\wt}{\widetilde}

         \newcommand{\pr}{^{\prime}}

         \newcommand{\ce}{\colonequals}
         \newcommand{\ov}{\overline}
         \newcommand{\un}{\underline}
         \newcommand{\sm}{\mathrm{sm}}
         
         \renewcommand{\b}{\textbf}
         \newcommand{\surjects}{\twoheadrightarrow}
         
         \newcommand{\isoto}{\overset{\sim}{\longrightarrow}}
         \newcommand{\isomfrom}{\overset{\sim}{\longleftarrow}}
         \newcommand{\et}{\mathrm{\acute{e}t}}	                                           % for etale cohomology (mainly used in subscripts)*/
         \newcommand{\alg}{\mathrm{alg}}		                                               % algebraic closure (mainly used in superscripts)*/
         \newcommand{\Div}{{\mathrm{Div}}}	                                               % Cartier divisor associated to a generically exact complex*/
         \newcommand{\aff}{^{\mathrm{aff}}}

         \providecommand{\fps}[1]{[\![#1]\!]}
         \providecommand{\lps}[1]{(\!(#1)\!)}
         \providecommand{\SP}[1]{\cite{SP}*{\href{https://stacks.math.columbia.edu/tag/#1}{#1}}}

         \providecommand{\f}[2]{\frac{#1}{#2}}
         \providecommand{\fps}[1]{\llbracket#1\rrbracket}
         \providecommand{\lps}[1]{(\!(#1)\!)}

\newextarrow{\xbigtoto}{{15}{15}{15}{12}}
   {\bigRelbar\bigRelbar{\bigtwoarrowsleft\rightarrow\rightarrow}}
         
         \DeclareMathOperator{\Ker}{Ker}			                       % Kernel*/
         \DeclareMathOperator{\Spec}{Spec}		                       % Spectrum of a ring*/
         \DeclareMathOperator{\Proj}{Proj}		                       % Proj of a graded ring*/
         \DeclareMathOperator{\Hom}{Hom}			                       % Set of arrows between two object*/
         \DeclareMathOperator{\Ass}{Ass}			                       % Associated primes*/
         \DeclareMathOperator{\Frac}{Frac}		                       % Field of fractions*/
         \DeclareMathOperator{\depth}{depth}		                       % Depth of a module*/
         \DeclareMathOperator{\Supp}{Supp}		                       % Support of a function*/
         \DeclareMathOperator{\id}{id}			                       % identity*/
         \DeclareMathOperator{\cube}{cube}
         \DeclareMathOperator{\sq}{sq}
         \DeclareMathOperator{\Ext}{Ext}		                           	% Derived functors of Hom*/
         \DeclareMathOperator{\Set}{\textbf{Set}}
         \DeclareMathOperator{\Norm}{Norm}		                                                  % Norm*/
         \DeclareMathOperator{\Pic}{Pic}		                                                  % Picard group*/
         \DeclareMathOperator{\Cl}{Cl}		                                                  % Class group*/
         \DeclareMathOperator{\codim}{codim}		                                                  % codimension*/

         \newcommand{\ba}{\begin{aligned}}
         \newcommand{\ea}{\end{aligned}}
         \newcommand{\be}{\begin{equation}}
         \newcommand{\ee}{\end{equation}}
         \newcommand{\pf}{\begin{proof}}
         \newcommand{\bpf}{\begin{proof}}
         \newcommand{\epf}{\end{proof}}
         \newcommand{\bsol}{\begin{solution}}
         \newcommand{\esol}{\end{solution}}
         \newcommand{\bthm}{\begin{thm}}
         \newcommand{\ethm}{\end{thm}}

         \newcommand{\bprop}{\begin{prop}}
         \newcommand{\eprop}{\end{prop}}
         \newcommand{\bcor}{\begin{cor}}
         \newcommand{\ecor}{\end{cor}}

         \newcommand{\brem}{\begin{rem}}
         \newcommand{\erem}{\end{rem}}

         \newcommand{\brems}{%
  \refstepcounter{thm}% 步进主定理计数器（例如占用 2.2 的位置）
  \begin{rems}\mbox{}% 调用已有的无编号 Remarks 环境，\mbox{} 用于换行
  \begin{enumerate}[label=\thethm.\arabic*, ref=\thethm.\arabic*]
}
\newcommand{\erems}{%
  \end{enumerate}
  \end{rems}
}

         \newcommand{\begs}{\begin{egs} \hfill \begin{enumerate}[label=\b{\thenumberingbase.},ref=\thenumberingbase]}
         
         \newcommand{\eegs}{\end{enumerate} \end{egs}}
         \newcommand{\eremstweak}{\end{enumerate} \end{rems-tweak}}
         \newcommand{\eremst}{\end{enumerate} \end{rems-tweak}}
         \newcommand{\blem}{\begin{lemma}}
         \newcommand{\elem}{\end{lemma}}

         \newcommand{\bconj}{\begin{conj}}
         \newcommand{\econj}{\end{conj}}
         \newcommand{\bprob}{\begin{Problem}}
         \newcommand{\eprob}{\end{Problem}}

         \newcommand{\bq}{\begin{Q}}
         \newcommand{\eq}{\end{Q}}
         \newcommand{\benum}{\begin{enumerate}[label={{\upshape(\alph*)}}]}
         \newcommand{\benuma}{\begin{enumerate}[label={{\upshape(\arabic*)}}]}
         \newcommand{\benumr}{\begin{enumerate}[label={{\upshape(\roman*)}}]}
         \newcommand{\eenum}{\end{enumerate}}
         \newcommand{\bc}{\begin{comment}}
         \newcommand{\ec}{\end{comment}}
         \newcommand{\bd}{\begin{defn}}
         \newcommand{\ed}{\end{defn}}
         \newcommand{\bque}{\begin{que}}
         \newcommand{\eque}{\end{que}}
         \newcommand{\bfct}{\begin{fact}}
         \newcommand{\efct}{\end{fact}}

         \newcommand{\beg}{\begin{eg}}
         \newcommand{\eeg}{\end{eg}}

         \newcommand{\bcl}{\begin{claim}}
         \newcommand{\ecl}{\end{claim}}

         \newcommand{\x}{\text}

         \newcommand{\q}{\quad}
         \newcommand{\qq}{\quad\quad}

         \newcommand{\tst}{\textstyle}

         \newcommand{\sdiv}{{\mathscr{D}\text{\!\!\!\:\:\raisebox{-.13ex}{\emph{iv}}}}}
         \newcommand{\xmapsfrom}[1]{%
    \mathrel{\reflectbox{$\xmapsto{\reflectbox{$\scriptstyle#1$}}$}}%
}

\tikzset{
    labl/.style={anchor=south, rotate=90, inner sep=.5mm}
}

\usepackage{aliascnt}
\newtheoremstyle{thms-style}
  {3.5pt}{3.5pt}      % 上下间距
  {\itshape}          % 正文字体 (斜体)
  {}                  % 缩进
  {\bfseries}         % 标题字体 (加粗)
  {.}                 % 标点
  { }                 % 标点后空格
  {}                  % 头部定义

\newtheoremstyle{remark-style}
  {3.5pt}{3.5pt}
  {\rmfamily}         % 正文字体 (罗马体，非斜体)
  {}
  {\itshape}          % 标题字体 (斜体，不加粗)
  {.}
  { }
  {}

\newtheoremstyle{claim-style}
  {3pt}{3pt}
  {\itshape}
  {}
  {\itshape}          % 标题字体 (斜体，这里也可以改成 \bfseries)
  {.}
  { }
  {}

\theoremstyle{thms-style}
\newtheorem{thm}{Theorem}[section] 
\Crefname{thm}{Theorem}{Theorems} % 定义引用名称

\newtheorem{mainthm}{Theorem}

\crefname{mainthm}{Theorem}{Theorems}

\newaliascnt{prop_cnt}{thm} 
\newtheorem{prop}[prop_cnt]{Proposition}
\aliascntresetthe{prop_cnt}
\Crefname{prop_cnt}{Proposition}{Propositions}

\newaliascnt{lemma_cnt}{thm}
\newtheorem{lemma}[lemma_cnt]{Lemma}
\aliascntresetthe{lemma_cnt}
\Crefname{lemma_cnt}{Lemma}{Lemmas}

\newaliascnt{cor_cnt}{thm}
\newtheorem{cor}[cor_cnt]{Corollary}
\aliascntresetthe{cor_cnt}
\Crefname{cor_cnt}{Corollary}{Corollaries}

\theoremstyle{remark-style} % 切换样式

\newaliascnt{rem_cnt}{thm}

\aliascntresetthe{rem_cnt}
\Crefname{rem_cnt}{Remark}{Remarks}

\newtheorem{rem}[rem_cnt]{Remark}  % 定义 rem，共享 rem_cnt 计数器
\Crefname{rem}{Remark}{Remarks}     % 确保引用 rem 时也显示 "Remark"

\newtheorem*{rems}{Remarks}

\theoremstyle{claim-style} % 切换样式

\newtheorem{claim}{Claim}[thm]
\Crefname{claim}{Claim}{Claims}

\newtheorem{cl-tweak}[thm]{Claim} % 注意：如果希望 cl-tweak 也是 3.2.1 格式，请改用下一行：
\theoremstyle{thms-style} % 切回主样式
\newtheorem{conj}[thm]{Conjecture} % 简单的共享，不处理别名(引用时会显示 Theorem)
\newtheorem{Q}[thm]{Question}
\newtheorem{variant}[thm]{Variant}

\newtheoremstyle{defs}{3.5pt}{3.5pt}{}{}{\bfseries}{.}{ }{}
\theoremstyle{defs}
\newtheorem{defn}[thm]{Definition}

\newtheorem{eg}[thm]{Example}

\newtheorem*{egs}{Examples}

\Crefname{claim}{Claim}{Claims}
\Crefname{sublemma}{Lemma}{Lemmas}
\Crefname{conj}{Conjecture}{Conjectures}
\Crefname{cor}{Corollary}{Corollaries}
\Crefname{defn}{Definition}{Definitions}
\Crefname{eg}{Example}{Examples}
\Crefname{prop}{Proposition}{Propositions}
\Crefname{Q}{Question}{Questions}
\Crefname{rem}{Remark}{Remarks}
\Crefname{thm}{Theorem}{Theorems}
\Crefname{variant}{Variant}{Variants}

\theoremstyle{thms}
\newtheorem{thm-tweak}[subsection]{Theorem}
\Crefname{thm-tweak}{Theorem}{Theorems}
\newtheorem{lemma-tweak}[subsection]{Lemma}
\Crefname{lemma-tweak}{Lemma}{Lemmas}
\newtheorem{cor-tweak}[subsection]{Corollary}
\Crefname{cor-tweak}{Corollary}{Corollaries}
\newtheorem{prop-tweak}[subsection]{Proposition}
\Crefname{prop-tweak}{Proposition}{Propositions}
\newtheorem{conj-tweak}[subsection]{Conjecture}
\Crefname{conj-tweak}{Conjecture}{Conjectures}

\theoremstyle{defs}
\newtheorem{defn-tweak}[subsection]{Definition}
\Crefname{defn-tweak}{Definition}{Definitions}
\newtheorem{eg-tweak}[subsection]{Example}
\Crefname{eg-tweak}{Example}{Examples}
\newtheorem*{rems-tweak}{Remarks}
\newtheorem{rem-tweak}[subsection]{Remark}
\Crefname{rem-tweak}{Remark}{Remarks}

\newtheoremstyle{subsection-tweak}
   {3.5pt}
   {3.5pt}%
   {}
   {}%
   {\bfseries}
   {}%
   {.5em}
   {\thmnumber{\@{#1}{}\@{#2}.}%
    \thmnote{~{\bfseries#3.}}}

\theoremstyle{subsection-tweak}
\newtheorem{pp}[thm]{}
\newcommand{\bpp}{\begin{pp}}
\newcommand{\epp}{\end{pp}}
\newtheorem{pp-t}[subsubsection]{}

\theoremstyle{subsection-tweak}

\theoremstyle{subsection-tweak}
\newtheorem{pp-tweak}{}

      \makeatletter
      \def\@tocline#1#2#3#4#5#6#7{
          \begingroup
          \@ifempty{#4}{}{}

          \parindent\z@ \leftskip#3\relax \advance\leftskip\@tempdima\relax
          #5\hskip-\@tempdima
            \ifcase #1
             \or\or \hskip 2em \or \hskip 1em \else \hskip 3em \fi%
            #6\nobreak\relax
          \dotfill\hbox to\@pnumwidth{\@tocpagenum{#7}}\par
          \nobreak
          \endgroup
        }
       \def\l@section{\@tocline{1}{0pt}{1pc}{}{}}

      \renewcommand{\tocsection}[3]{%
        \indentlabel{\@ifnotempty{#2}{\makebox[1.3em][l]{%
          \ignorespaces#1 \bfseries{#2}.\hfill}}}\bfseries{#3}
          \vspace{-3.5pt}}

      \renewcommand{\tocsubsection}[3]{%
        \indentlabel{\@ifnotempty{#2}{\hspace*{-0.5em}\makebox[2.1em][l]{%
          \ignorespaces#1#2.\hfill}}}#3
          \vspace{-4.5pt}}

\makeatother

\makeatletter % This is for proper formatting of appendices in the table of contents
\newcommand\appendix@section[1]{%
  \refstepcounter{section}%
  \orig@section*{Appendix \@Alph\c@section. #1}%
}
\let\orig@section\section
\g@addto@macro\appendix{\let\section\appendix@section}
\makeatother

      \DeclareMathVersion{normal2}

\author{Ning Guo}
\address{Institute for Advanced Study in Mathematics, Harbin Institute of Technology, Xidazhi 92, 150001 Harbin, China}
\email{ninguo@hit.edu.cn}
\date{\today}

\def\UTFviii@defined#1{%
  \ifx#1\relax
      \PackageError{inputenc}{Unicode\space char\space\expandafter
                              \UTFviii@splitcsname\string#1\relax
                              \MessageBreak
                              not\space set\space up\space
                              for\space use\space with\space LaTeX}\@eha
  \else\expandafter
    #1%
  \fi
}

\makeatletter
\def\UTFviii@defined#1{%
  \ifx#1\relax
      ?%
  \else\expandafter
    #1%
  \fi
}

\makeatother


\begin{document}

\title{Divisors on coherent schemes and homogeneous spaces}

\vspace{-25pt}

\begin{abstract}
We investigate the positivity and extension of invertible sheaves on group homogeneous spaces over coherent bases. Bypassing the failure of standard limit arguments and the classical Weil--Cartier correspondence, we develop a valuative divisor theory on locally coherent schemes. This establishes an exact correspondence between effective valuative divisors and rank-one reflexive sheaves, yielding a non-Noetherian Ramanujam--Samuel theorem. To homologically control special fibre degenerations, we study morphisms of (N)-type; these govern the descent of generically trivial invertible sheaves and establish the theorems of the cube and the square without smoothness hypotheses. Utilizing the Picard-admissibility of group actions, we construct ample invertible sheaves explicitly from one-codimensional orbit boundaries. This achieves the rigid extension of generic polarizations to integral models over Pr\"ufer bases, structurally generalizing Raynaud's classical proof of his quasi-projectivity theorems.
\end{abstract}

\maketitle

% \vspace{-40pt}
\hypersetup{
    linktoc=page,     %set to all if you want both sections and subsections linked
}

\renewcommand*\contentsname{}
\q\\
\tableofcontents

\section{Introduction}
%For a group scheme $G$ over a scheme $S$, an $S$-scheme $X$ has a $G$-action signifies a morphism of $S$-group functors $G\ra \un{\Aut}(X)$. 
%Such $X$ is an \emph{fppf $G$-homogeneous space} (resp., an \emph{fppf $G$-torsor}) if  further 
%\benumr
%\item the structural morphism $X\ra S$ is an fppf cover, and
%\item the morphism $G\times_SX\ra X\times_S X$ via $(g,x)\mapsto (g\cdot x, x)$ is an fppf cover (resp., an isomorphism).
%\eenum

Let $G$ be a group scheme over a scheme $S$. 
Recall that an fppf cover $X \ra S$ equipped with a $G$-action is a $G$-\emph{homogeneous space} (resp.\ a $G$-\emph{torsor}) if the canonical morphism $G\times_S X\ra X\times_S X$ is an fppf cover (resp.\ an isomorphism). 
In his book \cite{Ray70}, Raynaud systematically investigated the ampleness of invertible sheaves on homogeneous spaces, establishing the following quasi-projectivity theorem:
\bthm[\cite{Ray70}*{chapitre~VI, corollaire~2.5}]\label{Ray-main}
For a normal scheme $S$ and an $S$-smooth group scheme $G$ with connected $S$-fibres, every $G$-homogeneous space $X$ is $S$-locally quasi-projective\footnote{A morphism $f\colon X\ra S$ of schemes is \emph{locally quasi-projective} if it is quasi-projective locally on $S$.}.
\ethm

\Cref{Ray-main} extends the classical theorem of Chow \cite{Cho57} asserting the quasi-projectivity of connected smooth algebraic groups (the smoothness hypothesis later removed by Conrad \cite{CGP15}*{Proposition~A.3.5} via descent techniques).
To establish this, Raynaud explicitly constructed invertible sheaves from positive cycles of codimension one, subsequently applying his own ampleness criteria.
In the process, he obtained extension theorems for (semi)ampleness\footnote{An invertible sheaf $\sL$ is \emph{semiample} if there is an integer $n>0$ such that $\sL^{\otimes n}$ is generated by global sections.} on homogeneous spaces over Dedekind base schemes. 
%The aim of this article is to acquire analogues of aforementioned Raynaud's results over Pr\"ufer base schemes, namely, the schemes whose all local rings are valuation rings. 
%However, in such course of his approach to \Cref{Ray-main}, Raynaud repeatedly imposed Noetherian assumption to the base schemes, some of which have obstructions to be removed by using limit argument. 

While Raynaud's approach is effective in classical settings, arithmetic geometry increasingly demands a non-Noetherian foundation. For instance, the base rings in the perfectoid theory \cite{Sch12} and arc-topology \cite{BM21}  are valuation rings of higher ranks or nondiscrete, such as $\sO_{\bC_p}$ and $k\fps{u}+vk\lps{u}\fps{u}\subset k\lps{u}\lps{v}$; or Pr\"ufer rings, whose local rings are valuation rings, \emph{e.g.} the ring of algebraic integers $\ov{\bZ}$. 
Doing algebraic geometry over such bases runs into a fundamental obstacle. The classical Weil--Cartier correspondence breaks down. 
Furthermore, the standard limit arguments are notoriously unfeasible here; a limit of positive cycles of codimension $1$ does not, in general, correspond to a Cartier divisor on the  limit scheme.

To circumvent these obstructions, we propose that  \emph{locally coherent} schemes provide the natural setting where non-Noetherian pathologies are homologically controlled. 
A ring is \emph{coherent} if its finitely generated ideals are finitely presented; geometrically, the structure sheaf is coherent as a module over itself. 
Indeed, this notion has been increasingly recognized, for instance, in $p$-adic algebraic K-theory \cite{AMM22}.
As established in \cite{GL24}*{Theorem~2.21}, on any locally coherent, topologically locally Noetherian normal scheme, every local ring is either a valuation ring or has depth $\ge 2$. 
This theorem enables a dévissage: the global study of divisorial phenomena on coherent schemes reduces to local questions over valuation rings.

Building upon this dévissage, we replace the discrete $\mathbb{Z}$-multiplicities of classical Weil divisors with the totally ordered value groups $\Gamma_x$ of the governing valuation rings. We introduce the notion of \emph{valuative Weil divisors} and establish the following  correspondences that subsume the classical Weil--Cartier theory.

\begin{mainthm}
Let $X$ be a locally coherent, topologically locally Noetherian scheme.
\benumr
\item \textup{(\Cref{Weil-Rflx}\ref{Weil-Rflx-cor})} If $X$ is normal, there is an exact, order-preserving bijection between the monoid of effective valuative divisors $\mathfrak{Z}^1_+(X)$ and the monoid of rank-one reflexive ideal sheaves $\mathfrak{R}(X)$. 
\item \textup{(\Cref{regular-Weil-Cartier})} If $X$ is regular, then we have the Weil--Cartier correspondence $\operatorname{Cl}(X) \simeq \operatorname{Div}(X)$.
\eenum
\end{mainthm}

Empowered by this framework, we obtain a non-Noetherian counterpart to the Ramanujam--Samuel theorem (\Cref{Ramanujam--Samuel}): for smooth schemes over normal bases, relative effective Cartier divisors correspond precisely to flat one-codimensional cycles. 
Crucially, over Prüferian bases, this further ensures that effective Cartier divisors on generic fibres uniquely extend to the whole space via schematic closures. 

Applying this extension to the quasi-projectivity of homogeneous spaces, we bypass the Noetherianness  classically required by Raynaud to control orbit boundaries. In our coherent setting, without invoking finiteness conditions, this approach constructs ample invertible sheaves directly from the orbit geometry:

\begin{mainthm}[Effective ampleness via group orbit boundary, \Cref{constr-ample}]\label{mt1}
Let 
\begin{itemize}
\item $S$ be a locally coherent, topologically locally Noetherian, normal scheme.
\item $G$ be an $S$-flat group scheme of finite type with connected fibres,
\item $X$ an $S$-smooth scheme with $G$-action and an $S$-quasi-affine open $U$.
\end{itemize}
If $G_s$ is smooth whenever $\sO_{S,s}$ is a valuation ring, then one-codimensional irreducible components $D_i$ of the orbit boundary $G\cdot U-U$ span an effective divisor $D\ce \sum_{n_i\in \bZ_{>0}}n_iD_i$, and $\sL\ce \sO_X(D)$ is $S$-ample.
\end{mainthm}

Beyond the explicit generation of global positivity in \Cref{mt1}, a complementary issue is the \emph{persistence of polarizations}. When an integral model  is polarized over its generic point, extending this generic ampleness is severely obstructed by the failure of standard spreading-out arguments. Exploiting the schematic closures inherent to our valuative framework, we establish precise extension criteria.

\begin{mainthm}[Rigid extension of polarizations, \Cref{extend-extend}]\label{mt2}
Let $S$ be an affine integral Pr\"ufer scheme (e.g. $\Spec \ov{\bZ}$) with generic point $\eta$, and let $G$ be an $S$-smooth group  with connected fibres. Let $X$ be an $S$-separated, quasi-compact $G$-homogeneous space. For any invertible $\sO_{X_{\eta}}$-module $\sL_{\eta}$, the following hold:
\benumr
\item if $\sL_{\eta}$ is generated by global sections, then it extends to an $S$-semiample invertible sheaf on $X$;
\item if $\sL_{\eta}$ is ample, then there exists an integer $n\in \mathbb{Z}_{>0}$ such that $\sL_{\eta}^{\otimes n}$ extends to an $S$-ample invertible sheaf on $X$. 
Furthermore, any invertible extension of $\sL_{\eta}$ to $X$ is automatically $S$-ample.
\eenum
\end{mainthm}

While Theorems \ref{mt1} and \ref{mt2} recover Raynaud's results \cite{Ray70}*{chapitre~V, théorème~3.10 and chapitre~VIII, théorème~2} in the Noetherian setting, establishing these theorems over locally coherent bases strictly prohibits the standard limit arguments. The transition is achieved by resolving two geometric obstructions:

\begin{itemize}[leftmargin=*]
\item \emph{The divisorial limit obstruction.} The correspondence between flat pure codimension-one subschemes and relative effective Cartier divisors systematically fails to survive limit transitions. We bypass such obstruction, relying instead on our results of the valuative divisors and Ramanujam--Samuel theorem.
\item \emph{The fibre degeneration obstruction.} In arithmetic geometry, the degeneration of fibres is ubiquitous.
The inherent jumping of fibre properties renders standard limit arguments inadequate. 
Recognizing that such degenerations are unavoidable, we do not artificially restrict to smooth settings. Instead, we systematically develop morphisms of \emph{(N)-type} in Section~\ref{sec-N}, which emerge as the natural setting.
\end{itemize}

\bpp[(N)-type schemes]\label{def-N}
A flat, locally finitely presented scheme morphism $X\ra S$ is \emph{of (N)-type}, if each generic fibre $X_{\eta}$ is geometrically normal and any $X_s$ is geometrically reduced whenever $\sO_{S,s}$ is a valuation ring.
Such morphisms emerge for many moduli and local problems, as demonstrated below:
\benumr
\item Consider the universal curve $\pi\colon \overline{\mathcal{C}}_g \to \overline{\mathcal{M}}_g$ over the Deligne--Mumford moduli stack of stable curves of genus $g \ge 2$. 
The morphism $\pi$ is proper, flat, and representable by schemes. 
Because the generic family is smooth and the special fibres are stable curves possessing at worst ordinary nodes (which are geometrically reduced), every pulled-back scheme family is globally of (N)-type.
\item Schubert schemes $X_w \subset G/B$ over $\mathbb{Z}$ are  of (N)-type. As closures of Bruhat cells, they possess severe singularities along boundary strata, yet their geometric fibres remain normal and integral.
\item Semistable degenerations are of (N)-type, as special fibres are reduced normal crossing divisors.
\eenum
\epp

Over a coherent normal base, we prove that for suitable morphisms between (N)-type schemes, any generically trivial invertible sheaf on the source descends to the base (Corollaries \ref{N-pullback} and \ref{ext-for-N}). 
Consequently, as Example~\ref{eg-val} shows, for an (N)-type scheme $X$ with integral fibres over a valuation ring $V$, 
\[
\text{the restriction \q $H^1_{\et}(X,\mathbb{G}_m)\hra H^1_{\et}(X_{\Frac V}, \mathbb{G}_m)$\q is injective. }
\]
Being normal (\Cref{normality}) and closed under fibre products, schemes of (N)-type supply the nice category required to establish the theorems of the cube and the square in the absence of smooth hypotheses. 
When a group scheme acts on $X$, the interaction between the group actions and the Picard group is governed by the following notion (termed ``being of (C)-type'' by Raynaud).

\bpp[Picard-admissibility] For an $S$-scheme $X$ equipped with an action by an $S$-group scheme $G$, every invertible $\sO_X$-module $\sL$ induces a map of group functors by taking an absolute shift of $G$-translations
\[
\delta_{G,\sL}\colon G\ra \Pic_{X/S},\qq g\mapsto {}^{g}\sL\otimes \sL^{-1}.
\] 
The invertible sheaf $\sL$ is \emph{Picard-admissible} (for the $G$-action) if  $\delta_{G,\sL}$ is a homomorphism of groups.
\epp
By definition, Picard-admissibility is the functorial manifestation of the theorem of the square. Utilizing the aforementioned descent properties of (N)-type schemes, we systematically reduce the Picard-admissibility of an invertible sheaf to its generic fibre (\Cref{C-criterion}). 
This reduction is pivotal: it translates the group actions into criteria for (semi)ampleness (\Cref{homo-ample-crit}), which ultimately validate the boundary construction of ample sheaves in \Cref{mt1} and execute the rigid extension of polarizations in \Cref{mt2}.

%\bque
%By adding assumptions if necessary, can we generalize \Cref{triviality} to vector bundles, or, to $G$-torsors, where $G$ is an $S$-reductive group scheme?
%\eque

\bpp[Organization of the article]
Section~\ref{Hartogs} develops the valuative divisor theory on locally coherent schemes, establishing the correspondence between effective valuative divisors and rank-one reflexive sheaves.
Section~\ref{sec-ext} proves the non-Noetherian Ramanujam--Samuel theorem and the fpqc descent of generically trivial invertible sheaves.
Section~\ref{sec-N} introduces (N)-type morphisms to tame the degenerations of special fibres.
Section~\ref{sec-picard} leverages this to demonstrate the theorems of the cube and the square without smoothness, thereby reducing Picard-admissibility to generic fibres.
Section~\ref{sec-crit} translates this admissibility into geometric (semi)ampleness criteria, culminating in the boundary construction of ample sheaves (\Cref{mt1}).
Section~\ref{sec-ext-ample} concludes by the extension of generic polarizations over Pr\"ufer bases (\Cref{mt2}).
\epp
\bpp[Notations and conventions]
For a ring $R$, we denote its weak dimension by $\mathrm{wdim}(R)$. For a scheme $X$, let $X^{(i)}$ denote the set of points of codimension $i$.
\epp
\subsection*{Acknowledgements} I thank K\k{e}stutis \v{C}esnavi\v{c}ius for motivating me to write this article.
Besides, K\k{e}stutis \v{C}esnavi\v{c}ius, Fei Liu, Arnab Kundu, and Weizhe Zheng read an earlier version, I thank them for their reading and advices. 
This project has received funding from the European Research Council (ERC) under the European Union's Horizon 2020 research, the  innovation programme (grant No.~851146), the grant 075-15-2022-289, the National Natural Science Foundation of China (grant No.~12501016), and the excellent environment for research of the IASM of Harbin Institute of Technology.

% \blem\label{mb}
% For a valuation ring $V$, a decent algebraic space $X$ flat and locally of finite type over $V$, and a generic point $\eta$ of the closed fibre $X_s$, the Henselian local ring $\sO^h_{X,\eta}$ is a valuation ring $W$ such that
% \[
%    \text{$W$ dominates $V$\q and \q the value groups are equal $\GG_W=\GG_V$.} 
% \]
% \elem
% \bpf
% By construction of $\sO^h_{X,\eta}$, after reducing to the case of essentially finite type flat $V$-algebras, the assertion follows from a result of Moret-Bailly \cite{MB22}*{Thm.~A}.
% \epf

\section{Divisors on coherent schemes}\label{Hartogs}

\bpp[Cartier Divisors]
Given a ringed space $(X,\sO_X)$, let $\sM_X$ denote the sheaf of germs of meromorphic functions.
The \emph{sheaf of divisors} is the quotient $\sdiv_X\ce \sM_X^{\ast}/\sO_X^{\ast}$, where $\sM_X^{\ast}\subset \sM_X$ denotes the subsheaf of invertible functions.
The group  $\Div(X)\ce \GG(X,\sdiv_X)$ consists of \emph{(Cartier) divisors} on $X$.
The canonical map $\GG(X,\sM_X^{\ast}) \to \Div(X)$ sends a function $f$ to its \emph{principal divisor} $\mathrm{div}(f)$.
The \emph{support} of a divisor $D$ is the closed subset $\Supp(D)\ce \{x\in X\mid D_x\neq 0\}$.
We say a divisor $D$ is \emph{effective}, denoted by $D\geq 0$, if it lies in the image of the subsheaf of regular sections $\sO_X\cap \sM_X^{\ast} \to \sdiv_X$.
This induces a partial order on $\Div(X)$ such that for any $f\in \GG(X,\sM_X^{\ast})$, we have $\mathrm{div}(f)\geq 0$ if and only if $f\in \GG(X,\sO_X)$.
\epp
We introduce the notion of generalized Weil divisors or valuative Weil divisors, which naturally extend the classical theory by taking coefficients in the value groups of the structure sheaf.
\bpp[Valuative (Weil) divisors]
For an integral ringed space $(X,\sO_X)$, we define the set of \emph{valuative points} of $X$ as $X_{\mathrm{Val}}^{(1)}\ce \{x\in X \mid \sO_{X,x} \text{ is a valuation ring}\}$ and the group $\mathfrak{Z}^1(X)$ of \emph{valuative divisors} as
\[
\tst \mathfrak{Z}^1(X)\ce \bigoplus_{x\in X_{\mathrm{Val}}^{(1)}}\GG_x,
\]
where  $\GG_x$ is a totally ordered abelian group of the valuation ring $\sO_{X,x}$ with valuation $v_x: \sM_X(X)^\times \to \Gamma_x$.
A valuative divisor $D = (d_x)_{x} \in \mathfrak{Z}^1(X)$ is called \textit{effective}, denoted by $D \ge 0$, if $d_x \ge 0$ for all $x\in X^{(1)}_{\mathrm{Val}}$. 
The subgroup of effective valuative divisors is denoted by $\mathfrak{Z}^1_+(X)$.
If $X$ is topologically locally Noetherian, then the associated  \emph{principal divisor} of a nonzero $f \in \sM_X(X)^\times$ is defined as the following sum
\[
\tst \mathrm{val}(f) \coloneqq \sum_{x\in X^{(1)}_{\mathrm{Val}}} v_x(f) \cdot [x],
\]
where $[x]$ denotes the prime cycle of the closure $\ov{x}$.
The sum $\mathrm{val}(f)$ is well-defined owing to the topological assumptions on $X$.
Two valuative divisors $D_1, D_2$ are said to be \emph{linearly equivalent}, written $D_1 \sim D_2$, if their difference is a principal valuative divisor. 
We define the \emph{valuative divisor class group} as the quotient 
\[
\Cl(X) \coloneqq \mathfrak{Z}^1(X) / \mathrm{val}(\sM_X(X)^\times).
\]
For any finitely generated ideal sheaf $\cI\subset \sO_X$, its valuation at $x\in X^{(1)}_{\mathrm{Val}}$ is $v_x(\cI)\ce \min\{v_x(g)\mid g\in \cI_x\}$.
As finitely generated ideals in a valuation ring are principal,  $v_x(\cI)$ is well-defined. 
There is a map
\[
\tst \Psi\colon \cI\mapsto \mathrm{val}(\cI)\ce \sum_{x\in X^{(1)}_{\mathrm{Val}}}v_x(\cI)\cdot [x]
\]
that sends finitely generated ideal sheaves to valuative divisors. 
Conversely, every effective valuative divisor $D=(d_x)_x\in \mathfrak{Z}_+^1(X)$ is associated to a subsheaf $\sO_X(-D)\subset \sM_X$ as the following map shows
\[
\Phi\colon D\mapsto \sO_X(-D)\ce \{f\in \sM_X\mid  v_x(f)\geq d_x \text{ for all } x\in X^{(1)}_{\mathrm{Val}}\,\}.
\]
\epp

The maps $\Phi$ and $\Psi$ leverage coherence to resolve technical pathologies in the non-Noetherian world where maximal ideals are not principal.
Since finitely generated ideals in the valuation rings are principal, these objects are faithfully captured by the value groups $\GG_x$. This identifies reflexive rank-one sheaves as the natural generalization of classical Weil divisors, leading to the following correspondence theorem.

\bthm\label{Weil-Rflx} 
Let $X$ be a locally coherent, topologically locally Noetherian, normal scheme.
\benumr
\item For every finitely generated ideal $\cI\subset \sO_X$, its composed image $\Phi(\Psi(\cI))=\cI^{\vee\vee}$ is the reflexive hull.
\item\label{Weil-Rflx-cor} Let $\mathfrak{R}(X)$ denote the monoid of reflexive ideal sheaves on $X$. 
There is an order-preserving bijection
\begin{equation*}
    % --- 左侧：集合对应 ---
    \begin{tikzcd}[ampersand replacement=\&, column sep=2.5em]
        \mathfrak{R}(X) \arrow[r, shift left=0.8ex, "\Psi"] 
        \& \mathfrak{Z}^1_+(X) \arrow[l, shift left=0.8ex, "\Phi"]
    \end{tikzcd}
    \qquad \qquad 
    % --- 右侧：具体构造（单行显示） ---
    \cI \xmapsto{\Psi} \operatorname{val}(\cI)
    \;, \qquad
    \sO_X(-D) \xmapsfrom{\Phi} D
\end{equation*}
between monoids $\mathfrak{R}(X)$ and $\mathfrak{Z}^1_+(X)$.
Furthermore, this correspondence descends to an isomorphism between $\operatorname{Cl}(X)$ and the group $\sO_X\x{-}\mathbf{Rflx}^1$ of isomorphism classes of rank-one reflexive sheaves.
\eenum
\ethm
\bpf 
Since $v_x(\cI)$ is well-defined, so is the map $\Psi \colon \mathfrak{R}(X)\to \mathfrak{Z}^1_+(X)$.
Given a divisor $D=(d_x)_x\in \mathfrak{Z}^1_+(X)$, since $d_x\geq 0$, the image $\Phi(D)$ is an ideal sheaf.
The image $\Phi(D)=\bigcap_{x\in X^{(1)}_{\mathrm{Val}}}\{x\in \sM_X\mid v_x(f)\geq d_x\}$.
Since $D \in \mathfrak{Z}^1_+(X)$ belongs to the direct sum, $d_x > 0$ for only locally finitely many valuative points. Hence, locally, $\Phi(D)$ is a finite intersection of $\sO_X$ and finitely many coherent submodules of $\sM_X$. Thus, being a finite intersection of coherent modules, $\Phi(D)$ is finitely generated and coherent.
Then \cite{GL24}*{Proposition~2.22} yields the reflexivity of $\Phi(D)$, so the map $\Phi\colon \mathfrak{Z}^1_+(X)\to \mathfrak{R}(X)$ is well-defined.

For a reflexive ideal sheaf $\cI$, we show that $\Ass(\sO_X/\cI)$ consists of valuative points. 
Pick $\eta\in \Ass(\sO_X/\cI)$ and let $R\ce \sO_{X,\eta}$.
Suppose that $R$ is not a valuation ring, then \cite{GL24}*{Theorem~2.21} implies that $\depth(R) \ge 2$. 
The quotient $Q \ce R/\cI_\eta$ is supported at $\fm_R$.
    Consider the short exact sequence $0 \to \cI_\eta \to R \to Q \to 0$.
    Applying the functor $\Hom_R(\cdot, R)$ to this sequence  yields the long exact sequence:
    \[
        0 \to \Hom_R(Q, R) \to \Hom_R(R, R) \to \Hom_R(\cI_\eta, R) \to \Ext^1_R(Q, R) \to \cdots
    \]
    We have $\Hom_R(R/\cI_{\eta}, R) =\Ext^1_R(Q, R) = 0$. 
    Therefore, the natural map $\Hom_R(R, R) \isoto  \Hom_R(\cI_\eta, R)$ is an isomorphism.
  Since $\cI_\eta$ is a reflexive $R$-module, we have $\cI_\eta \cong ((\cI_\eta)^\vee)^\vee \cong R^\vee \cong R$.
    Thus, $\cI_\eta=(f)$ is principal for a nonzerodivisor $f \in \fm_R$ such that $\sqrt{(f)} = \fm_R$. 
    In particular, every element in $\fm_R$ is nilpotent modulo $f$, so $\depth_R(R/f)=0$.
    Consequently, the depth formula \cite{GL24}*{Lemma~2.6} gives $\depth_R(R)=1+\depth_R(R/f)=1$, which contradicts the assumption that $\depth_R R\geq 2$. 

Let $\xi_i\in \Ass(\sO_X/\cI)$.
By construction, we write the effective valuative divisor of $\cI$ as a sum
\[
\tst \Psi(\cI)=D_{\cI}\ce \operatorname{val}(\cI)=\sum_i v_{\xi_i}(\cI)\cdot [\xi_i].
\]

We verify that $\Phi$ and $\Psi$ are inverse to each other.
First, consider the composite $\Phi \circ \Psi$. 
By \cite{GL24}*{Proposition~2.22}, we have $\cI=\bigcap_{x\in X^{(1)}_{\mathrm{Val}}}\cI_{x}$.
Since $\cI$ is trivial beyond $V(\cI)$, this intersection only occurs at those $\xi_i\in \Ass(\sO_X/\cI)$.
Note that $\cI_{\xi_i}$ is principal and isomorphic to $\sO_{X,\xi_i}$ by multiplying a generator of valuation $v_{\xi_i}(\cI)$, hence the construction $\Phi(D_{\cI})$ just yields $\cI$ itself. 
This shows that $\Phi(\Psi(\cI))=\cI$.

Conversely, consider $\Psi \circ \Phi$. 
Let $D = \sum d_x [x] \in \mathfrak{Z}^1_+(X)$ and let $\mathcal{J} \ce \sO_X(-D) = \Phi(D)$.
We need to check that $v_x(\mathcal{J}) = d_x$ for all $x \in X^{(1)}_{\mathrm{val}}$.
Fix a valuative point $x_0\in X^{(1)}_{\mathrm{Val}}$. The stalk of $\mathcal{J}$ at $x_0$ is given by
\[
   \tst \mathcal{J}_{x_0} = \left( \bigcap_{y \in X^{(1)}_{\mathrm{val}}} \{ f \in \sM_X \mid v_y(f) \ge d_y \} \right)_{x_0}.
\]
Since specializations in $X^{(1)}_{\mathrm{Val}}$ correspond to quotients of valuation groups (\cite{Guo24}*{Proposition~A.2(v)}), we are reduced to considering independent valuation rings.
In such case, the condition $v_{x_0}(f) \ge d_{x_0}$ is the only constraint local to $x_0$, we have $\mathcal{J}_{x_0} = \{ f \in \sO_{X,x_0} \mid v_{x_0}(f) \ge d_{x_0} \}$, whose valuation in  $\sO_{X,x_0}$  is exactly $d_{x_0}$.
Thus $\Psi(\Phi(D)) = D$.
This establishes the desired bijection between $\mathfrak{R}(X)$ and $\mathfrak{Z}^1_+(X)$.

Finally, we extend this to the divisor class group.
The bijection extends to a group homomorphism by defining $\operatorname{val}(\cI)$ for a fractional reflexive sheaf $\cI$ via linearity.
The monoid structure on $\mathfrak{Z}^1_+(X)$ is addition $D_1 + D_2$ and the corresponding structure on $\mathfrak{R}(X)$ is the reflexive product $(\cI \otimes \mathcal{J})^{\vee\vee}$.
For principal divisors, if $f \in K(X)^\times$, then $\sO_X(-\operatorname{val}(f))$ is the fractional ideal generated by $f$, which is isomorphic to $\sO_X$.
Conversely, if $\cI \cong \sO_X$ is a free rank-one sheaf, then it is generated by a unit section, implying its associated divisor is trivial in $\operatorname{Cl}(X)$.
Thus, the correspondence descends to an isomorphism
\[
\operatorname{Cl}(X) \cong \sO_X\x{-}\mathbf{Rflx}^1. \qedhere
\]
\epf 
We generalize the classical result that Weil divisors are Cartier divisors on regular schemes as follows.

\bthm\label{regular-Weil-Cartier} 
On any locally coherent, regular, topologically locally Noetherian scheme $X$, rank-one reflexive sheaves are exactly invertible sheaves.
In particular, we have the following isomorphisms
\[
\x{$\sO_X\x{-}\mathbf{Rflx}^1\simeq \Cl(X)\simeq \Div(X)$\q and \q $\Cl_+(X)\simeq \Div_+(X)$.}
\]
\ethm
\bpf 
For a rank-one reflexive sheaf $\sF$ on $X$, its invertible locus $U$ is nonempty and open. 
As the rank function is constructible, the open immersion $j\colon U\hra X$ is quasicompact.
By \cite{GL24}*{Theorem~2.19}, the complement $X\backslash U$ has weak dimension $\geq 2$, hence \emph{op.~cit.}~Theorem~2.20 and its proof yield 
\[
   \sF\isoto j_{\ast}j^{\ast}\sF\isoto j_{\ast}\det(j^{\ast}\sF[0])\isoto j_{\ast}j^{\ast}\det \sF[0]\isomfrom \det \sF[0],
\]
where $\sF[0]$ is a perfect complex thanks to the coherent regularity of $X$.
Hence $\sF$ is invertible.
\epf

\section{Extension and descent of divisors}\label{sec-ext}
Unless otherwise specified, all divisors in this section refer to Cartier divisors.
We study extension and descent of divisors and invertible sheaves on normal schemes.
Specifically, taking schematic closure of effective divisors on generic fibres yields relative effective divisors on Pr\"ufer bases as the unique flat model.

\blem[\cite{RG71}*{corollaire~3.4.7}]\label{raygru} 
Every algebra  flat of finite type over a domain is finitely presented.
\elem

\blem\label{flat-equiv}
Let $f\colon X\ra S$ be a morphism locally of finite type where $S$ is integral.
For the conditions (i)~$X\ra S$ is flat;
(ii)~$\sO_X$ is $\sO_S$-torsionfree;
(iii)~$X$ is the schematic closure of its generic fibre $X_{\eta}$; and 
(iv)~$\Ass(\sO_X)\subset X_{\eta}$,
we have (i) $\Rightarrow$ (ii) $\Leftrightarrow$ (iii) $\Leftrightarrow$ (iv).
If $S$ is Pr\"ufer, then all conditions are equivalent.
\elem

\begin{proof}[Proof of \Cref{flat-equiv}] 
When $S$ is Pr\"ufer, we have (ii) $\Rightarrow$ (i).
In general cases, the implications (i) $\Rightarrow$ (ii) $\Leftrightarrow$ (iii) are standard. 
For (iii) $\Leftrightarrow$ (iv), recall \SP{05C3} that $a \in \sO_S$ is a zerodivisor on $\sO_X$ if and only if it belongs to $\Ass(\sO_X)$.
Consequently, (iii)  is equivalent to the condition that for every $\mathfrak{p} \in \text{Ass}(\sO_X)$, the contraction $\mathfrak{p} \cap \sO_S$ is the zero ideal (0).
This means $f(\mathfrak{p}) = \eta$, which is precisely condition (iv).
\end{proof}

%%%%Now, we recall the following criterion for parafactoriality in a relative case.
%%%%\blem[\cite{EGAIV4}*{proposition~21.14.3}]\label{RS-classical}
%%%%For a normal scheme $S$ and an $S$-smooth scheme $X$,
%%%%\benumr
%%%%\item\label{RS-cl} if $S$ is locally Noetherian, then any one-codimensional cycle $Z$ on $X$ such that $\Supp(Z)$ avoids all maximal points of $S$-fibres of $X$ is locally principal;
%%%%\item\label{GR-ext} for a closed subset $Y\subset X$ satisfying the following condition (\ref{GR}), the pair $(X,Y)$ is parafactorial:
%%%%\begin{equation}\label{GR}
%%%%\text{$\codim(Y_{\eta},X_{\eta})\geq 2$ for every generic point $\eta\in S$\q  and \q $\codim(Y_s, X_s)\geq 1$ for every $s\in S$.}\tag{PF}
%%%%\end{equation}
%%%%\eenum
%%%%\elem

\bprop\label{Ramanujam--Samuel} 
For a normal domain $R$, a smooth $R$-scheme $X$ of finite presentation, and a closed subscheme $Z\subset X$, consider the following conditions:
\benumr 
\item\label{flat-codim-one} $Z$ is $R$-flat, finitely presented, and of pure codimension $1$;
\item\label{rel-eff-cart} $Z$ is a relative effective divisor on $X$;
\item\label{aff-large} $Z$ has no embedded associated points such that $j\colon X\backslash Z\hra X$ is affine and $R$-fibrewise dense;
\item\label{schcl-effcart} $Z$ is the schematic closure in $X$ of an effective divisor $Z_{\eta}\subset X_{\eta}$.
\eenum
Then \ref{flat-codim-one}$\Lra$\ref{rel-eff-cart}$\Lra$\ref{aff-large}$\Rightarrow$\ref{schcl-effcart}.
If we assume that $R$ is Pr\"uferian or $Z$ is $R$-flat, then \ref{flat-codim-one}--\ref{schcl-effcart} are equivalent. 
\eprop

\bpf 
The implication \ref{rel-eff-cart} $\Rightarrow$ \ref{schcl-effcart} follows from \Cref{flat-equiv}.
Then, we show that \ref{schcl-effcart}$\Rightarrow$\ref{rel-eff-cart} when $R$ is Pr\"uferian. 
In fact, instead of considering smooth morphisms, we only need the regularity of $X$.

\bcl\label{effcart-schcl-releffcart} 
For a locally coherent regular scheme $X$ flat of finite type over an integral Pr\"ufer scheme $S$, the schematic closure $D$ of an effective divisor $D_{\eta}$ on the generic fibre $X_{\eta}$ is a relative effective divisor.
\ecl
\begin{proof}[Proof of the claim]

Since $S$ is Prüferian, the schematic closure $D$ of $D_\eta$ is $S$-torsion-free, hence $S$-flat (by \Cref{flat-equiv}). 
Because $D$ is $S$-flat, its formation commutes with flat base change, so we may localize to assume that $S=\Spec R$ for  a valuation ring $R$.
For the same reason, as $R$ is a filtered direct union of valuation subrings $R_i$ of finite ranks such that  $R_i\hra R$ are faithfully flat, we may assume that $\dim R<\infty$.

We may assume that $X$ is connected.
The rank of $R$ is finite, so $X$ is a topologically locally Noetherian.
Let $\eta\in X^{(1)}_{\text{Val}}$ be a valuative point, then the stalk $(\cI_D)_x$ of the ideal sheaf $\cI_D$ is determined by a quotient map $V\ra V/a$, where $V$ is a valuation ring and $a\in \fm_V$.
Since $D_{\eta}$ is an effective divisor, the ideal sheaf $\cI_D$ is indeed reflexive by \cite{GL24}*{Proposition~2.22}.
In particular, we have $\cI_D\in \sO_X\x{-}\mathbf{Rflx}^1$.
The coherent regularity of $X$ and \Cref{regular-Weil-Cartier} imply that $\cI_D$ is invertible such that $D$ is an effective divisor on $X$.
Besides, the $R$-flatness of $D$ follows from \Cref{flat-equiv}.
Consequently, $D$ is a relative effective divisor. \qedhere
\end{proof}

When $Z$ is $R$-flat, the implication \ref{schcl-effcart}$\Rightarrow$\ref{rel-eff-cart} follows from a similar argument for \ref{flat-codim-one}$\Rightarrow$\ref{rel-eff-cart} below: we base change everything over a valuation ring to show that the special fibre is an effective divisor.

Write $R$ as a filtered direct limit of Noetherian normal subdomains $R_{\lambda}$ (with affine transitions), we may descend $X$ and $Z$ to $X_{0}$ and $Z_{0}$ over some $R_{0}$. 
It remains to descend  \ref{flat-codim-one}--\ref{aff-large} to the Noetherian base. 
For \ref{flat-codim-one}, we apply \cite{EGAIV3}*{théorème~11.2.6} for the $R$-flatness. 
Using \cite{EGAIV3}*{proposition~14.3.13 and théorème~13.1.3}, each nonempty fibre $Z_s$ is of pure codimension one. 
So this fibrewise codimension condition descends by \cite{EGAIV2}*{corollaire~6.1.4}.   
After descent, the closed $Z_0$ is of pure codimension one in $X_0$.
For \ref{rel-eff-cart}, the invertibility of the ideal sheaf $\cI_Z$ descends by \cite{EGAIV3}*{théorème~8.5.2 and proposition~8.5.5}.
The descent for \ref{aff-large} follows from \cite{EGAIV3}*{théorème~8.3.11} for constructible closed subschemes ($X$ is quasi-separated so $j$ is retrocompact). 
Now, we are reduced to the Noetherian case.

\ref{flat-codim-one}$\Rightarrow$\ref{rel-eff-cart}.
By using the fibral criterion \SP{062Y}, it suffices to show that each fibre $Z_s\neq \emptyset$ is an effective divisor on $X_s$. 
Let $V$ be a (discrete) valuation ring dominating $R$ centered at $s$.
Since the residue field of $V$ is just $\kappa(s)$, after base changing we may assume that $R=V$.
It is clear that the generic fiber $Z_{\eta}$ is an effective divisor.
Hence the schematic closure $Z$ of $Z_{\eta}$  by \ref{schcl-effcart}$\Rightarrow$\ref{rel-eff-cart} is a relative effective divisor.

\ref{rel-eff-cart}$\Rightarrow$\ref{aff-large}.
If \ref{rel-eff-cart} holds, then $Z$ has no embedded associated points and the open immersion $j$ is affine (\SP{07ZU}). 
Besides, each fibre $Z_s$ is of codimension $\geq 1$ in $X_s$. 
Consequently, we get \ref{rel-eff-cart}$\Rightarrow$\ref{aff-large}.

\ref{aff-large}$\Rightarrow$\ref{flat-codim-one}.
By \cite{EGAIV4}*{corollaire~21.12.7}, the closed subscheme $Z$ is of pure codimension $1$.
As the immersion $X\backslash Z\hra X$ is dense, the cycle $Z$ is locally principal due to  Ramanujam--Samuel \cite{EGAIV4}*{proposition~21.14.3}, hence is an effective divisor. 
The $R$-flatness of $Z$ then follows from \SP{046Z}.
\epf

\brems
\item If $X$ has regular fibres instead of being $R$-smooth, then the implications above are valid except \ref{aff-large}$\Rightarrow$\ref{flat-codim-one}: though $Z$ is a Weil divisor in the Noetherian case, its being locally principal need that all nonvaluative points of $X$ are parafactorial.  
However, when $R$ is Pr\"uferian, we still have the parafactoriality of local rings at $x\in X$ when $\dim \sO_{X,x}\geq 2$, see \cite{GL24}*{Theorem~3.11}. 
\item For non-Pr\"uferian $R$, \ref{schcl-effcart}$\Rightarrow$\ref{rel-eff-cart} is false without the $R$-flatness of $Z$.
For example, when $R=k[x,y]$ with fraction field $K$ and $X=\Spec k[x,y,t]$, the schematic closure $Z$ of $Z_{\eta}\ce K[t]/(xt-y)$ is not $R$-flat and the fibre of $Z$ at $(0,0)$ is $\bA^1_k$, so $Z$ is not a relative effective divisor. 
\erems

% \bcor 
% For an affine open immersion $j\colon U\hra X$, every irreducible component of $X\backslash U$ is either a relatively effective divisor on $X$, or the closure $\ov{\{\xi\}}$ of a maximal point $\xi$ of $S$-fibre of $X$.
% \ecor
% \bpf 
% Note that $X$ is topologically Noetherian, by \SP{01OX}, every open subset of $X$ is retrocompact, in particular, $j$ is quasi-compact.
% If $U$ contains all the maximal points of $S$-fibres of $X$, then \Cref{Ramanujam--Samuel} suffices.
% Assume that there is a maximal point $\eta$ of $S$-fibres of $X$ contained in $Z\ce X\backslash U$. 
% By \Cref{mb}, the points of $X$ that specialize to $\xi$ are also maximal points of $S$-fibres of $S$. 
% Hence there is a unique irreducible component $Z_{\alpha}\subset Z$ with generic point $\xi_{\alpha}$ containing $\eta$ such that $Z_{\alpha}=\ov{\{\xi_{\alpha}\}}$.
% Let $U\pr\ce X\backslash \ov{(Z\backslash Z_{\alpha})}$, then $j\pr\colon U\hra U\pr$ is a quasi-compact affine open immersion and we are reduced to prove the assertion for $U\pr\backslash U$. Since there are finitely many maximal points of $S$-fibres of $X$, we reduce iteratively to the resolved case \Cref{Ramanujam--Samuel}.
% \epf

\blem\label{detect-positive} 
Let $X$ be a  scheme. For every divisor $D$ on $X$ and the points $x\in X$, we have 
\[
\text{$D\geq 0$ (resp., $D=0$) \q iff \q $D_x\geq 0$ (resp., $D_x=0$) whenever $\depth \sO_{X,x}=1$.}
\]
\elem
\bpf 
We prove that, $D\geq 0$ follows from $D_x\geq 0$ for all $x\in X$ such that $\GG(\Spec \sO_{X,x}, \sO_X)\not\simeq \GG(U_x,\sO_X)$, where $U_x$ is the punctured spectrum.
The effectiveness of $D$ is a local property, so we assume that $D=\mathrm{div}(f)$ for a regular meromorphic function $f$ on $X$. 
The effective locus $\mathrm{dom}(f)\ce \{x\in X\,|\, \mathrm{div}(f)_x\geq 0\}$ is open.
Hence, by hypothesis, for $T\ce X\backslash \mathrm{dom}(f)$, we have an isomorphism $\GG(X,\sO_X)\simeq \GG(X\backslash T, \sO_X)$.
Consequently, there is an $\wt{f}\in \GG(X,\sO_X)$ such that $\wt{f}|_{\mathrm{dom}(f)}=f$, which forces $T=\emptyset$, so $D$ is effective.
\epf

% \beg
% For a valuation ring $V$ with spectrum $S$, we may choose $\mathbf{\Omega}(S)\ce \{s\in S\,| \,\dim \sO_{S,s}\geq 1\}$.
% For a locally Noetherian normal scheme $X$, we may choose $\mathbf{\Omega}(X)\ce \{x\in X\,|\, \dim \sO_{X,x}=1\}$.
% \eeg
%\bd
%For a normal scheme $S$ and an $S$-flat finitely presented scheme $X$ such that 
%\benumr
%\item for every generic point $\eta\in S$, the fibre $X_{\eta}$ is geometrically normal, and
%\item for every point $s\in S$, the fibre $X_s$ is geometrically reduced,
%\eenum
%
%\ed
Recall that a faithfully flat morphism of schemes induces an injective pullback map on divisors, hence we may view the image of a pullback as a divisor subgroup and form the quotient divisor group.
\bprop\label{inv-div} 
For an fpqc morphism $f\colon X\ra Y$ of schemes, a divisor on $X$ descends to $Y$ if and only if it does so around every $y\in Y$ with $\mathrm{depth}(\sO_{Y,y})\leq 1$. 
More precisely, the following map
\[
   \text{$\frac{\Div(X)}{f^{\ast}\Div(Y)}\hra \prod_{\x{$\mathrm{depth}(\sO_{Y,y})\leq 1$}} \frac{\Div(X(y))}{f^{\ast}_y\Div(Y(y))}$}
\]
is injective, where we denote $Y(y)\ce \Spec \sO_{Y,y}$, and $X(y)\ce X\times_{Y}Y(y)$.

\eprop
\bpf 
For a divisor $D$ on $X$, we need to show that it is $f^{\ast}(\Delta)$ for a divisor $\Delta$ on $Y$ if and only if it does so over all $y\in Y$ of depth $\leq 1$.
The necessity follows from the transitivity of taking preimages of divisors (\cite{EGAIV4}*{21.4.4}).
It remains to show the sufficiency.
For an open subset $U\subset Y$, by faithful flatness, taking inverse image along $f_U$ induces an injective map $\Div(U)\hra \Div(f^{-1}(U))$.
Therefore, for two open subsets $U_1$ and $U_2$ of $Y$, if there are divisors $\Delta_1$ on $U_1$ and $\Delta_2$ on $U_2$ such that $D|_{f^{-1}(U_i)}=f_{U_i}^{\ast}(\Delta_i)$ for $i=1, 2$, then $\Delta_1$ and $\Delta_2$ coincide over $U_1\cap U_2$ so that the solution extends.
Suppose that there is a largest open $U \subsetneq Y$ such that $D|_{f^{-1}(U)}=f^{\ast}_{U}(\Delta|_{U})$.
It suffices to prove that any maximal point $y\in Y\backslash U$ has a neighbourhood over which $D$ descends, thereby violating the maximality of $U$.

As $f$ is quasi-compact, a limit argument \cite{EGAIV4}*{proposition~20.3.8(ii)} reduces us to showing that $D$ descends to $\Delta$ over $\Spec \sO_{Y,y}$.
By hypothesis, we may assume that $Y=\Spec \sO_{Y,y}$ has depth $\geq 2$ at $y$ and $U=Y\backslash \{y\}$. 
The open immersion $j_X\colon X_{U}\hra X$ is the base change of $j_Y\colon U\hra Y$ and $j_X\colon X_{U}\hra X$ along $f$.
The divisor $\Delta|_{U}$ corresponds to a pair $(\sL_{U},s_U)$, where $\sL_{U}$ is an invertible $\sO_{U}$-module and $s_U$ is a regular meromorphic section of $\sL_{U}$ on $U$.
Similarly, we write $D=(\sL_X,s_X)$.
The goal is to descend  $(\sL_X,s_X)$ over $Y$.
By hypothesis, we have $j_X^{\ast}(\sL_X)=f^{\ast}_{U}(\sL_{U})$, hence taking $(j_X)_{\ast}$ yields $(j_X)_{\ast}j_X^{\ast}\sL_X=(j_X)_{\ast}f_{U}^{\ast}(\sL_{U})$.
Note that $f$ is flat, by depth formula \cite{GR18}*{Corollary~10.4.46} we have $\mathrm{depth}(\sO_{X,x})\geq 2$ for every $x\in f^{-1}(y)$,  so flat base change yields $\sL_X=(j_X)_{\ast}f_{U}^{\ast}(\sL_{U})=f^{\ast}(j_Y)_{\ast}(\sL_{U})$. 
Since $Y$ is local, $\sL_X$ descends to the trivial line bundle $(j_Y)_{\ast}(\sL_{U})\simeq \sO_Y$.

It remains to descend $s_X$ to a regular meromorphic section of $\sO_Y$ on $Y$.
Since $U\subset Y$ is schematically dense, by \cite{EGAIV4}*{proposition~20.2.11}, the section $s_U$ extends to a regular meromorphic section $s_Y$ such that $s_Y|_{U}\circ f_U=s_X|_{X_{U}}$.
In conclusion, the divisor $\Delta\ce (\sO_Y,s_Y)$ satisfies $f^{\ast}(\Delta)=D$.
\epf
\bcor\label{div-gen-zero} 
Let $f\colon X\ra Y$ be an fpqc morphism of locally coherent schemes.
If $Y$ is normal and for every $y\in Y$ with $\operatorname{wdim}(\sO_{Y,y})=1$ the fibre $X_y$ is integral, then the following sequence is exact:
\[
   \x{$0\ra \Div(Y)\overset{f^{\ast}}{\ra} \Div(X)\ra \Div(X)|_{Y^{(0)}}\ra 0$.}
\]
Namely, a divisor $D$ on $X$ is a pullback from $Y$ if and only if $D_{\xi}=0$ for every maximal point $\xi\in Y$.
\ecor

Reducing to $y\in Y$ with $\sO_{Y,y}$ a nontrivial valuation ring, Moret-Bailly's theorem yields a weakly unramified extension at the generic point of $X_y$, hence $\GG_{\sO_{Y,y}}\simeq \GG_{\sO_{X_y,\zeta}}$, which is the key input for descending divisors.

\blem[\cite{MB22}*{th\'eor\`eme~A}]\label{mb}
Let $S$ be an irreducible Pr\"ufer scheme with generic point $\eta$, and let $X\to S$ be flat and of finite type.
For $s\in S$ and a generic point $\xi$ of a generically reduced irreducible component of $X_s$, the inclusion $\sO_{S,s}\subset \sO_{X,\xi}$ is weakly unramified, that is, we have \ $\GG_{\sO_{S,s}}\xrightarrow{\sim}\GG_{\sO_{X,\xi}}$.
\elem

\begin{proof}[Proof of \Cref{div-gen-zero}]
By \Cref{inv-div}, the problem is local on $y\in Y$ for which $\sO_{Y,y}$ is a valuation ring.
When $\sO_{Y,y}$ is a field, it is clear that $D_y$ is a pullback. 
Now assume that $\sO_{Y,y}$ is a nontrivial valuation ring. 
Let $\zeta\in X_y$ be the generic point.
By Moret-Bailly's \Cref{mb},  $\sO_{Y,y}\subset \sO_{X_y,\zeta}$ is a weakly unramified extension of valuation rings. 
Recall that a weakly unramified extension of valuation rings $V_1\subset V_2$ induces a bijection of value groups $\GG_{V_1}\simeq \GG_{V_2}$.
Every divisor $D_2$ on $\Spec V_2$ is a section of $\GG(\Spec V_2, \sM^{\ast}/\sO^{\ast})$, which, by Hilbert's 90, just equals to $\GG_{V_2}$, so $D_2$ is a pullback from $D_1\in \GG_{V_1}$.
\end{proof}

\bcor\label{gen-trivial} 
Under the setup of \Cref{div-gen-zero}, for maximal points $\eta\in Y$, we have
\[
 \tst  H^1(Y,\bG_m)=\Ker(H^1(X,\bG_m)\to \prod_{\eta}H^1(X_{\eta},\bG_m)),
\]
i.e., for an invertible $\sO_X$-module $\sL$, we have $\sL_{\eta}\simeq \sO_{X_{\eta}}$ iff it is a pullback of an invertible $\sO_Y$-module.
\ecor
\bpf 
It suffices to find a regular meromorphic section $s$ of $\sL$ on $X$ to form a divisor $D_X$ such that $(D_X)_{\eta}=0$, that is, $(D_X)_x=0$ for every point $x\in X_{\eta}$.
If so, then \Cref{div-gen-zero} yields a divisor $D_Y=(\sM, s\pr)$ such that $f^{\ast}D_Y\simeq D_X$ including an isomorphism $f^{\ast}\sM\simeq \sL$.

First, we reduce to the case when $Y$ is affine and integral.
By \cite{GR18}*{Corollary~10.5.9}, all points in $\Ass(\sO_X)$ are lying over $\eta$ so there is a dense open affine neighborhood $Y_0\subset Y$ of $\eta$ such that $f^{-1}(Y_0)\subset X$ is schematically dense thanks to \cite{GR18}*{Proposition~10.5.6(i)}.
Replacing $Y$ by $Y_0$ suffices.

Then, we find the section $s$.
Since $X_{\eta}$ is integral, $\Ass(\sO_{X})$ consists of a single point $\zeta$, the generic point of $X$.
For the isomorphism $\sL_{\eta}\isoto \sO_{X_{\eta}}$, let $s_{\eta}$ be the preimage of the unit section.
Let $X_0\subset X$ be an affine open neighborhood of $x\in X_{\eta}$ in $X$ and let $s$ be a section of $\sL$ on $X_0$ such that $s_{x}=(s_{\eta})_{x}$.
By \emph{op.~cit.}, the open subsets in $X$ that are schematically dense are exactly those that contains the set $\Ass(\sO_X)$, hence $X_0$ is schematically dense in $X$.
Thus, $s$ is a regular meromorphic section of $\sL$ on $X$.
Finally, by \cite{EGAIV4}*{propositions~20.2.11 and 20.3.5}, $s$ induces the section $s_{\eta}$ of $\sL_{\eta}$, as desired.
\epf

\beg[\emph{cf.} \SP{0EX8}]\label{eg-val}
For a valuation ring $V$ and a connected $V$-flat finite type scheme $X$, if all $V$-fibres of $X$ are reduced and the closed fibre is integral, then we have an injection
\[
   H^1(X,\mathbb{G}_m)\hra H^1(X_{\eta},\mathbb{G}_m).
\]
In fact, we may assume that $X\ra \Spec V$ is surjective. 
The connectedness of $X$ and Moret-Bailly's \Cref{mb} imply that all $Y$-fibres of $X$ are integral. 
By \Cref{gen-trivial}, every line bundle that trivializes on generic fibres is a pullback, which is trivial since $V$ is local.
\eeg

% \bprop\label{gen-tri-tri}
% For a valuation ring $V$ with spectrum $(S,\eta)$, a flat finite type morphism $f\colon X\ra S$ of irreducible schemes, and an $\sO_X$-invertible module $\sL$ on $X$,  assume that one of the following  holds
% \benumr
% \item\label{sm-int-fib} $X$ has integral geometrically reduced nongeneric fibres and geometrically normal generic fibre;
% \item\label{proper-sp} $f$ is proper with integral closed fibre such that $\GG(X_{\eta},\sO_{\eta})\simeq \Frac V$;
% \item\label{proper-int-fib} $f$ is proper with reduced fibres and integral closed fibre, and $V$ has finite rank.
% \eenum
% Then $\sL$ is trivial if and only if $\sL_{\eta}$ is trivial.
% In particular,  we have an injective map of pointed sets
% \[
% H^1_{\et}(X,\mathbb{G}_m)\hra H^1_{\et}(X_{\eta}, \mathbb{G}_m).
% \]
% \eprop
% \bpf
% The assertion \ref{sm-int-fib} just follows from \Cref{gen-trivial} when $Y=S$.
% The assertion \ref{proper-sp} is \SP{0EX8}.
% The assertion \ref{proper-int-fib} is essentially a geometric variant of \ref{proper-sp}. First, we prove the following claim

% Since $f$ has reduced fibres and $X$ is irreducible, by \Cref{discon-degen-nonred-discon}, the generic fibre $X_{\eta}$ is integral.
% Because $X_{\eta}$ is proper over $K=\Frac V$, the global section $\GG(X_{\eta},\sO_{X_{\eta}})$ is a finite-dimensional $K$-algebra, which, due to the integrality of $X_{\eta}$, is the reduced local Artin $K$-algebra $K$ itself.
% Then case \ref{proper-sp} applies.
% \epf
\bprop\label{general-div-seq} 
Let $f\colon X\ra S$ be a faithfully flat and finitely presented morphism  with geometrically integral fibres.
If $S$ is normal, quasi-compact and quasi-separated, then the following sequence is exact
\[
   \x{$0\ra \Div(S)\overset{f^{\ast}}{\ra} \Div(X)\ra \Div(X)|_{S^{(0)}}\ra 0$.}
\]
\eprop
\bpf 
It is clear that the image of $f^{\ast}$ vanishes on fibres at maximal points of $S$.
Suppose that a divisor $D$ on $X$ becomes trivial on $X_{\eta}$ for each maximal point $\eta\in S$.
The goal is to show that $D$ is a pullback.

As $S$ is quasi-compact, we may assume that $S$ is normal and integral. 
By a limit argument combining \SP{01ZA}, normalization of Nagata schemes (see \Cref{normality}\ref{weak-N-normal}) and gluing schemes, we write $S$ as a limit of normal integral Noetherian schemes $S_i$ with affine dominant transitions $t_{ji}\colon S_j\to S_i$.
Then $f\colon X\ra S$ descends to a faithfully flat finite type morphism $f_i\colon X_i\to S_i$ and the divisor $D$ descends to a divisor $D_i$ on $X_i$. 
Since the set of $s\in S_i$ such that $(X_i)_a$ is geometrically integral is locally constructible (\cite{EGAIV2}*{théorème~9.7.7(iv)}) and contains the image of $S$,  by \SP{05F4} we may assume that $f_i$ also has geometrically integral fibres. 
Next, we reduce to the case when $(D_i)_{\eta_i}=0$ for the generic point $\eta_i\in S_i$.
For every $j \geq i$, let $\eta_j$ be the unique generic point of $S_j$. 

\bcl 
Let $F_j \ce \{\eta_j\}$ if $(D_j)_{\eta_j} \neq 0$ and let $F_j \ce \emptyset$ otherwise.
Then, the sets with transitions $\{F_j,t_{kj}|_{F_k}\}_{j\geq i}$ form an inverse system of sets. 
\ecl
\bpf[Proof of the claim] 
Every transition $t_{kj}$ is affine dominant for $k\geq j\geq i$, hence induces a map $S_k^{(0)}\to S_j^{(0)}$.
We show that $t_{kj}|_{F_k}\colon F_k\to F_j$ is well-defined. 
Let $p_{kj}\colon X_k\to X_j$ be the base change of $t_{kj}$.
If $\eta_k\in F_k$ then $(D_k)_{\eta_k}\neq 0$, so its image $\eta_j\colon t_{kj}(\eta_k)\in S_j^{(0)}$ satisfies that $(D_k)_{\eta_k}=(p_{kj}|_{\eta_{k}})^{\ast}(D_j)_{\eta_j}$.
As $(p_{kj}|_{\eta_k})^{\ast}$ is induced by a field extension, it is injective on divisors. 
Therefore, if $\eta_j\not\in F_j$, then $(D_j)_{\eta_j}=0$ so that $(D_k)_{\eta_k}=0$, which leads to a contradiction with $\eta_k\in F_k$. 
So $\{F_j,t_{kj}|_{F_k}\}_{j\geq i}$ form an inverse system. 
\epf
Since $t_{ji}$ are dominant and affine, there is an isomorphism $X^{(0)}\isoto \varprojlim_{j\geq i}X_j^{(0)}$. 
In particular, by the claim, the inverse limit $F\ce \varprojlim_{j\geq i}F_j$ is the generic point $\eta\in S$ or $\emptyset$. 
However, by assumption $D|_{X_{\eta}}=0$ so the divisor trivializes over some $\eta_k$ for $k\geq i$, which contradicts the construction of $F_k$.
Hence, we have $F=\emptyset$ so that we may assume that $F_j=\emptyset$ for all $j\geq i$, or equivalently, that $(D_i)_{\eta_i}=0$ for the generic point $\eta_i\in S_i$.
 The problem is reduced to the Noetherian case solved in \Cref{div-gen-zero}.
\epf 
\bcor\label{general-ker-gen} 
Under the setup of \Cref{general-div-seq}, for maximal points $\eta\in S$, we have
\[
 \tst  H^1(S,\bG_m)=\Ker(H^1(X,\bG_m)\to \prod_{\eta}H^1(X_{\eta},\bG_m)),
\]
i.e., for an invertible $\sO_X$-module $\sL$, we have $\sL_{\eta}\simeq \sO_{X_{\eta}}$ iff it is a pullback of an invertible $\sO_S$-module.
\ecor

\section{Morphisms of (N)-type}\label{sec-N}

\bd 
Let $f\colon X\to S$ be a flat, locally of finite type morphism of schemes.
When $X$ and $S$ are locally coherent, we say that $f$ is of \emph{(N)-type} if its generic fibres are geometrically normal and the fibre $X_s$ is geometrically reduced whenever $\mathrm{wdim}\, \sO_{S,s}=1$.
For general schemes, we say that $f$ is of \emph{strict (N)-type} if all non-generic fibres are geometrically reduced and generic fibres are geometrically normal.
\ed
It is clear that all morphisms of strict (N)-type are of (N)-type.
We study their properties as follows.

% If $S$ is Noetherian, then the conditions of fibres are can be dropped over $s\in S$ of $\dim \sO_{S,s}\geq 2$.
% With the following properties, we study the theorems of cube and of square for (N)-type schemes.

\blem\label{normality} 
Let $f\colon X\ra S$ be a flat, locally of finite type morphism.
\benumr 
\item\label{N-normal} If $f$ is of (N)-type, $X$ is locally topologically Noetherian and $S$ is normal, then $X$ is normal.
\item\label{weak-N-normal} If $f$ is a strict (N)-type morphism of schemes, then $X$ is normal\footnote{If $S$ is locally coherent and normal (\emph{e.g.}, Noetherian normal), then for $X$ to be normal, the condition that the $S$-fibres $X_s$ are geometrically reduced is only required for $s\in S$ such that $\operatorname{wdim}(\sO_{S,s})=1$, see \cite{GL24}*{Theorem~2.21}. 
Raynaud's definition is restricted to the Noetherian case and fails to merge with our notion via limit arguments. Roughly speaking, (N)-type morphisms require that the locus of non-geometrically-reduced fibres has a ``sufficiently large'' weak codimension.}.
\item\label{N-fibre-product} (Strict) (N)-type morphisms are stable under flat base change: if $S\pr\to S$ is flat and $f$ is of (strict) (N)-type, then so is $X_{S\pr}\to S\pr$.
Moreover, being (strict) (N)-type is preserved under fibre products. 
\item\label{N-integral} In the case \ref{N-normal} or \ref{weak-N-normal}, assume further that $S$ is integral,  $X_1$ and $X_2$ have geometrically integral generic fibres, then $X_1\times_S X_2$ is integral and normal.
\eenum
\elem
\bpf
The assertion \ref{N-normal} follows from a coherent Serre's criterion for normality \cite{GL24}*{Theorem~2.21}.

For \ref{weak-N-normal}, we assume that $X$ and $S=\Spec R$ are affine connected.   
By \cite{EGAIV2}*{scholie~7.8.3(ii)(vi)} and a limit argument, we write $R$ as a filtered direct union of $R_0$-finite type normal subdomains $(R_{\lambda})_{\lambda\in \GGL}$ for a $\b{Z}$-finite type normal subdomain $R_0\subset R$ where $0\in \Lambda$ such that $X$ descends to an affine $R_{0}$-flat finite type scheme $X_{0}$.
The subset $E\subset \Spec(R_{0})$ over which $X_{0}\ra \Spec(R_{0})$ has geometrically reduced fibres is locally constructible (\SP{0579}) and contains the image of $\Spec R$ (\SP{0576}), hence by \SP{05F4} the images of $\Spec R_{\mu}\ra \Spec R_{0}$ are contained in $E$ for all large enough $\mu\in \GGL$. 
So $X$ descends to $X_{\mu}$, which is $R_{\mu}$-flat, of finite type, and all fibres are geometrically reduced.   
Since geometrically normality is preserved and reflected by base field extensions \SP{038P}, $X_{\mu}$ also has geometrically normal generic fibre. 
Then, we are reduced to the Noetherian case, where we combine Serre's criterion and the depth formula \cite{EGAIV2}*{proposition~6.3.1} for flat morphisms to check points in each $R_{\lambda}$-fibre, that for $\lambda\pr\geq \mu$ each $X_{\lambda\pr}$ satisfies $(R_1)+(S_2)$.
We conclude by taking the direct limit of normal domains \cite{EGAIV2}*{5.13.6}.

For \ref{N-fibre-product}, note that generic points of $S\pr$ lie over those of $S$, hence $X_{S\pr}$ has geometrically normal generic fibres.
By \cite{GL24}*{Lemma~A.3(ii)}, each $s\pr\in S\pr$ with $\mathrm{wdim}\,\sO_{S\pr,s\pr}=1$ has image $s$ such that $\mathrm{wdim}\, \sO_{S,s}\leq 1$, hence every $X_{s\pr}$ is geometrically reduced.
In particular, being (N)-type is stable under fibre products.

The assertion \ref{N-integral} follows from \cite{EGAIV2}*{proposition~4.6.21 and corollaire~2.3.5}: the fibre product $X_1\times_SX_2$ has geometrically integral generic fibre, hence is irreducible.
\epf
Given an integral scheme $S$ and a scheme morphism $f\colon X\ra S$, let $\rho(X)$ or $\rho_f(X)$ be the \emph{geometric number of irreducible components} of $X_{\eta}$, that is, the number of irreducible components of $X_{\eta}\otimes_{k_{\eta}}k_{\eta}^{\mathrm{alg}}$.
\blem\label{reduce-geonum-irrcomp} 
Let $f\colon X\ra S$ and $h\colon Y\ra S$ be (strict) (N)-type morphism of normal integral schemes.
Let $g\colon X\ra Y$ be an $S$-morphism and  $Z_0\subset X\times_S Y$ the connected component\footnote{The existence of $Z_0$ is guaranteed because $X\cong \GG_g\hra X\times_S Y$ is an immersion and $X$ is connected. } containing the graph $\GG_g$.
\benumr
\item\label{delta-one} For the composite $\pi\ce \mathrm{pr}_1|_{Z_{0}}\colon Z_{0}\hra X\times_S Y\surjects X$, we have $\rho_{\pi}(Z_{0})=1$.
\item\label{reduce-r} Let $X\times_S Y= Z_{0}\sqcup_i Z_i$ be the component decomposition. 
For the map $\pi_i\ce \mathrm{pr}_1|_{Z_i}$, we have 
$$\rho_{\pi_i}(Z_i)\leq \rho_{h}(Y)-1.$$
\eenum 
\elem
\bpf[Proof of the claim] 
Note that $Z_{0}$ is irreducible since $X\times_S Y$ is normal (\Cref{normality}\ref{N-normal} and \ref{weak-N-normal}).
Let $\eta$, $\xi$ and $\zeta$ be the generic points of $S$, $X$ and $Y$ respectively.
The scheme $X\times_S Y$ is $X$-flat, thus  its irreducible components correspond to those of $\xi\times_{\eta}\zeta$. 
In particular, $(Z_0)_{\xi}$ is connected, which, with the presence of a $\xi$-rational point and \cite{EGAIV2}*{corollaire~4.5.14}, yields $\rho_{\pi}(Z_0)=1$.
It remains to deduce \ref{reduce-r}.
In fact, $\rho_{\mathrm{pr}_1}(X\times_S Y)$ is equal to the number of irreducible components of $\Spec(k(\zeta)\otimes_{k(\eta)}k(\xi)^{\alg})$, so equals to $\rho_h(Y)$. 
Consequently, for any irreducible component $X_i$ of $(X\times_S Y)\backslash Z_0$, we get the desired inequality
\[
   \rho_{\pi_i}(X_i)\leq \rho_{\mathrm{pr}_1}(X\times_S Y)-\rho_{\pi}(Z_0)=\rho_h(Y)-1.\qedhere
\]
\epf

\blem\label{obstruction}
For a normal, locally coherent scheme $S$ and  a morphism $f\colon X\ra S$ of (N)-type, we have 
\[
   X^{(1)}_{\mathrm{Val}}=  \bigsqcup_{\eta\in S^{(0)}}X_{\eta}^{(1)}\cup \bigsqcup_{\mathrm{wdim}\sO_{S,s}=1}X_s^{(0)}.
\]
If $X$ is locally topologically Noetherian, then $x\in X^{(1)}_{\mathrm{Val}}$ iff $\mathrm{depth}\,\sO_{X,x}=1$, iff $\sO_{X,x}$ is a valuation ring.
\elem
\bpf 
We write $\mathfrak{X}_1\ce \bigsqcup_{\eta\in S^{(0)}}X_{\eta}^{(1)}$ and $\mathfrak{X}_2\ce \bigsqcup_{\mathrm{wdim}\sO_{S,s}=1}X_s^{(0)}$. 
Each local ring at $x\in \mathfrak{X}_1$ is a discrete valuation ring. 
If $x\in \mathfrak{X}_2$, then by Moret-Bailly's \Cref{mb}, $\sO_{X,x}$ is a valuation ring.
Conversely, every local ring $\sO_{X,x}$ satisfies the inequality of weak dimension \cite{GL24}*{Lemma~A.3}, hence $X^{(1)}_{\mathrm{Val}}\subset \mathfrak{X}_1\sqcup \mathfrak{X}_2$.

For the next assertion, by \Cref{normality}, $X$ is normal and locally coherent. 
So \cite{GL24}*{Theorem~2.21} (which needs topological Noetherianness) implies that $\mathrm{depth}\,\sO_{X,x}=1$ which amounts to that $\sO_{X,x}$ is a valuation ring.
Now assume that $\mathrm{depth}\,\sO_{X,x}=1$. 
If $f(x)=\eta$, then the normality of $X_{\eta}$ forces $x\in X_{\eta}^{(1)}$.
If $f(x)=s$ such that $\mathrm{wdim}\,\sO_{S,s}=1$, then the depth formula \cite{GR18}*{Theorem~10.4.37} forces that $x\in X_s^{(0)}$. 
It remains the case when $\mathrm{wdim}\, \sO_{S,f(x)}\geq 2$. 
Again by \cite{GL24}*{Theorem~2.21} we have $\mathrm{depth}\,\sO_{S,f(x)}\geq 2$ which contradicts with the depth formula. 
Consequently, $\mathrm{depth}\,\sO_{X,x}=1$ implies that $x\in X^{(1)}_{\mathrm{Val}}$.
\epf
In the view of \Cref{obstruction}, we get the following special cases of \Cref{detect-positive} and \Cref{div-gen-zero}.
\bcor\label{N-pullback}
For a morphism $f\colon X\ra S$ of (N)-type satisfying all conditions in \Cref{obstruction}, then
\[
\tst \text{$D\geq 0$ (resp., $D=0$) \q iff \q $D_x\geq 0$ (resp., $D_x=0$) for all $x\in X^{(1)}_{\mathrm{Val}}$.}
\]
\ecor

\bcor\label{ext-for-N} 
For a locally coherent normal scheme $S$, let $X$ and $Y$ be topologically locally Noetherian (N)-type $S$-schemes.
If for each $y\in Y^{(1)}_{\mathrm{Val}}$ the fibre $X_y$ is integral, then the following sequence is exact:
\[
   \x{$0\ra \Div(Y)\overset{f^{\ast}}{\ra} \Div(X)\ra \Div(X)|_{Y^{(0)}}\ra 0$.}
\]
In particular, for maximal points $\eta\in Y$, we have
\[
 \tst  H^1(Y,\bG_m)=\Ker(H^1(X,\bG_m)\to \prod_{\eta}H^1(X_{\eta},\bG_m)),
\]
i.e., for an invertible $\sO_X$-module $\sL$, we have $\sL_{\eta}\simeq \sO_{X_{\eta}}$ iff it is a pullback of an invertible $\sO_Y$-module.
\ecor
This exact sequence allows us to trivialize torsors over valuation rings, as the following example illustrates:
\beg 
If $V$ is a valuation ring and $X$ is a connected (N)-type $V$-scheme with integral closed fibre, then every line bundle on $X$ that trivializes over $X_{\Frac V}$ is trivial.
\eeg

Now we generalize the comparison \cite{GR18}*{Propositions~11.4.6} to (N)-type schemes over regular bases.
\bprop\label{positive-detect} 
Let $f\colon X\ra S$  an (N)-type morphism for regular, topologically locally Noetherian $S$.
For the smooth locus $X^{\mathrm{sm}}$ of $f$, taking restriction induces equivalences of categories
\[
   \x{ $\{\x{rank-one reflexive sheaves on } X\} \isoto \mathbf{Pic}(X^{\mathrm{sm}})$   \q and \q $\sO_X\x{-}\mathbf{Rflx}^1 \isoto \Div(X^{\mathrm{sm}})$.}
\]
\eprop
\bpf 
Suppose that $\sE$ is rank-one reflexive, then \Cref{regular-Weil-Cartier} implies that  $\sE|_{X^{\sm}}\in \mathbf{Pic}(X^{\sm})$, which guarantees the restriction functor is well-defined.
Each fibre $X_s$ is geometrically reduced, so is generically smooth. 
Similarly, generic fibres $X_{\eta}$ are smooth in codimension $1$.
We conclude that any local ring $\sO_{X,z}$ at $z\in X\backslash X^{\sm}$ is not a valuation ring, i.e., $\mathrm{wdim}(\sO_{X,z})\geq 2$.
Note that the open immersion $j\colon X^{\sm}\hra X$ is quasi-compact, by \cite{GL24}*{Theorem~2.20} the assertion holds. \qedhere

\epf

\section{Picard-admissibility of invertible sheaves with group actions}\label{sec-picard}
Leveraging the descent properties of (N)-type schemes, we establish the theorems of the cube and the square without smoothness hypotheses. This reduces the Picard-admissibility of an invertible sheaf to its generic fibre, preparing for global ampleness constructions.
\bpp[Theorems of cube and square]\label{thm-sq}
Let $S$ be a scheme and $f^i\colon X^i\to S$ ($i=1,2,3$) be $S$-schemes.
For an invertible sheaf $\sL$ on $C\ce X^1\times_S X^2\times_S X^3$, we introduce the \emph{cubical difference}
on $C\times_S C$ as follows. Identify $C\times_S C$ as the $S$-product $\prod_{i=1}^3X^i\times_S X^i$ and for each $I\subset\{1,2,3\}$ define the morphism 
\[
q_I\colon C\times_S C\ra C\q  \text{by \q $q_I|_{X^i\times_S X^i}=\Bigl\{$ \parbox{15em}{$\mathrm{pr}_1\colon X^i\times_SX^i\ra X^i$,\q if $i\in I$ \\  $\mathrm{pr}_2\colon X^i\times_S X^i\ra X^i$,\q if $i\not\in I$.}
}
\]
which prescribes the $i$-th factor.
Set
\[
\tst c(\sL)\;\ce\;\bigotimes_{I\subset\{1,2,3\}} q_I^{*}(\sL)^{\otimes (-1)^{|I|}}
\;\in\;\Pic(C\times_S C).
\]
We say that $\sL$ \emph{satisfies the theorem of cube}, or $\cube(\sL)$ holds, if $c(\sL)\simeq \sO_{C\times_S C}$.

Similarly, we put $W\ce X^1\times_S X^2$ and for $J\subset \{1,2\}$ define the following morphism 
\[
r_J\colon W\times_S W\times_S X^3\ra C\q  \text{by \q $r_J|_{X^i\times_S X^i}=\Bigl\{$ \parbox{15em}{$\mathrm{pr}_1\colon X^i\times_SX^i\ra X^i$,\q if $i\in J$ \\  $\mathrm{pr}_2\colon X^i\times_SX^i\ra X^i$,\q if $i\not\in J$}
}
\]
where the $X^3$-factor is the identity.
Let $\pi\colon W\times_S W\times_S X^3\to W\times_S W$ be the projection and set
\[
\tst s(\sL;X^1,X^2)\;\ce\;\bigotimes_{J\subset\{1,2\}} r_J^{*}(\sL)^{\otimes (-1)^{|J|}}
\;\in\;\Pic(W\times_S W\times_S X^3).
\]
We say that $\sL$ \emph{satisfies the theorem of square with respect to $(X^1,X^2)$},
and write $\sq(\sL;X^1,X^2)$, if 
\[
s(\sL;X^1,X^2)\;\simeq\;\pi^{*}\sM
\quad\text{for some }\sM\in\Pic(W\times_S W).
\]
The two notions are related by a formal identity: 
Indeed, writing $C\times_S C \simeq W\times_S W\times_S X^3\times_S X^3$,
let $\mathrm{pr}_k\colon C\times_S C\to W\times_S W\times_S X^3$ ($k=1,2$) be projections. Then one has a canonical isomorphism
\[
c(\sL)\;\simeq\;\mathrm{pr}_2^{*}s(\sL;X^1,X^2)\otimes \mathrm{pr}_1^{*}s(\sL;X^1,X^2)^{-1},
\]
hence $\sq(\sL;X^1,X^2)\Rightarrow \cube(\sL)$.
Conversely, If $\mathrm{cube}(\sL)$ holds and each of the following is satisfied
\benumr
\item $f^3\colon X^3\ra S$ permits a section;
\item $f^3$ is fpqc and for every $T\ra S$ we have $f^3_{T\ast}(\sO_{X^3_T})=\sO_T$\footnote{This is a stronger condition than that $f^3$ is cohomologically flat in degree zero.},
\eenum
then $\mathrm{sq}(\sL, X^1,X^2)$ also holds.  
\epp

The following generalizes \cite{BLR90}*{\S6.3, Theorem~1} and \cite{Ray70}*{chapitre~IV, théorème~2.3}.
\bprop\label{power-cube-square} 
Let $S$ be a normal integral scheme with generic point $\eta$.
Let $X^i\ra S$ be three scheme morphisms of strict (N)-type (resp., (N)-type).
Assume that $X^i_{\eta}$ and $X^3_s$ are geometrically integral (resp., when $\mathrm{wdim}\, \sO_{S,s}=1$), and $X^3\surjects S$ is surjective.
Let $\sL$ be an invertible sheaf  on $C=X^1\times_S X^2\times_S X^3$.
\benumr
\item\label{square-generic} For $W$ and $\pi$ as in \S\ref{thm-sq}, and the generic point $\xi\in W\times_S W$, we have the following equivalences 
\[
\text{$\mathrm{sq}(\sL,X^1,X^2)$ holds\q $\Leftrightarrow$ \q $\mathrm{sq}(\sL_{\eta}, X^1_{\eta}, X^2_{\eta})$ holds \q $\Leftrightarrow$ \q $s(\sL, X^1, X^2)_{\xi}$ is trivial.}
\]
\item\label{power-square} There is  $n\in \mathbb{Z}_{>0}$ such that $\mathrm{sq}(\sL^{\otimes n}, X^1,X^2)$ holds, a fortiori,  the theorem of cube holds for $\sL^{\otimes n}$.
\eenum
\eprop
\bpf 
For \ref{square-generic}, by definition, if $\mathrm{sq}(\sL,X^1,X^2)$ holds, then $\mathrm{sq}(\sL_{\eta}, X^1_{\eta}, X^2_{\eta})$ holds, and so $s(\sL,X^1,X^2)_{\xi}$ is a pullback of an invertible sheaf over $\xi$ hence trivial.
It remains to show that the triviality of $\mathrm{s}(\sL,X^1,X^2)_{\xi}$ implies that the theorem of square of $\sL$ for the factors $X^1$ and $X^2$, namely, $\mathrm{s}(\sL,X^1,X^2)$ is a pullback of an invertible sheaf along $\pi$.
Hence, we check the conditions in Corollaries~\ref{gen-trivial} and \ref{general-ker-gen}.
By \Cref{normality}\ref{N-fibre-product}, the fibre product $W\times_S W$ is of (strict) (N)-type over $S$.
Note that $\pi$ is fpqc.

For $w\in W\times_SW$ such that $\mathrm{wdim}\, \sO_{W\times_SW,w}=1$, by \cite{GL24}*{Lemma~A.3(ii)}, the image of $w$ in $S$ has local ring of weak dimension at most one, so $\pi^{-1}(w)$ is integral.
Consequently, Corollaries~\ref{gen-trivial} and \ref{general-ker-gen} apply and yield an invertible sheaf $\sM$ over $W\times_S W$ such that $\sL\simeq \pi^{\ast}(\sM)$, so $\mathrm{sq}(\sL,X^1,X^2)$ holds.

To prove \ref{power-square}, it suffices to exploit \ref{square-generic} to reduce us to the case when $S=\Spec k$ for a field $k$ and $(X^i)_{i=1}^3$ are locally of finite type, geometrically (normal and integral) $k$-schemes. 
For this, we first prove a claim.
\vskip 0.5cm
\bcl 
There is a finite field extension $K/k$ such that $\mathrm{sq}(\sL_K, X^1_K,X^2_K)$ holds.
\ecl 
By this claim, there is an invertible $\sO_{(W\times_k W)_K}$-module $\sM$ such that $s(\sL_K)\isoto \pi_K^{\ast}(\sM)$ is an isomorphism. 
So we have $s(\sL)^{\otimes r}=s(\sL^{\otimes r})\simeq \pi^{\ast}(\Norm_{K/k}(\sM))$, which leads to the assertion \ref{power-square}.
\epf

\bpf[Proof of the Claim] 
The assertion \ref{square-generic} permits us to replace $X^1$ by an affine open subset $U^1$. 
Up to a finite field extension, there is a $k$-rational point on $X^3$, hence we are reduced to showing that $\mathrm{cube}(\sL_{K}|_{U^1_K\times X^2_K\times X^3_K})$ holds for a finite extension $K/k$. 
By symmetry, it suffices to verify over $X^2_K\times X^3_K\times U^1_K$, or further, verify $\mathrm{sq}(\sL_{K}|_{X^2_K\times X^3_K\times U^1_K},X^2_K,X^3_K)$.
So we may assume that $X^3$ is affine.

Consider the perfect closure $k\subset k^{\mathrm{perf}}$ as a direct limit (union) of a finite subextensions, then a limit argument combined with faithfully flat descent reduce us to the case when $k$ is perfect.

As $X^3$ is separated and of finite type, Chow's lemma \cite{EGAII}*{5.6.1} yields  a projective birational map $g\colon \wt{X^3}\to X^3$.
By the normality of $X^3$, the map $g$ is an isomorphism on an open $U$ that contains all points of codimension $\leq 1$. 
The perfectness of $k$ ensures that the smooth and regular loci coincide, so we have $\codim(X^3\backslash U^{\mathrm{sm}})\geq 2$. 
Then Hartogs' extension permits us to replace $X^3$ by $U^{\mathrm{sm}}$, which is smooth and quasi-projective. 
Using \ref{square-generic} to shrink $X^1$ and $X^2$ to affine open subsets, we may assume that they are smooth and there are projective normal schemes $(Y^i)_{i=1}^3$ containing $(X^i)_{i=1}^3$ as open subschemes. 

Let $Z = X^3 \times_k X^3$, and consider the line bundle $\sN = p_{123}^{\ast}\sL \otimes p_{124}^{\ast}\sL^{-1}$ on $X^1 \times_k X^2 \times_k Z$, where $p_{123}$ and $p_{124}$ are the projections to the respective factors of $X^1 \times_k X^2 \times_k X^3 \times_k X^3_K$. Up to a finite field extension, we may choose $k$-rational points $x_1 \in X^1$ and $x_2 \in X^2$. 
Rigidifying $\sL$ along these points ensures that $\sN$ is trivialized along $\{x_1\} \times_k X^2 \times_k Z$ and $X^1 \times_k \{x_2\} \times_k Z$.
Let $E \subset Z$ be the maximal subscheme over which the rigidified line bundle $\sN$ is trivial. 
Because $(X^i)_{i=1}^2$ are dense opens of the projective normal schemes $(Y^i)_{i=1}^2$, standard representability properties of the Picard functor for projective normal schemes (combined with the seesaw principle \cite{CF23}*{Proposition~3.1}) guarantee that  $E$ is clopen in $Z$.
By construction, the pullback of $\sN$ along the diagonal $\mathrm{id}_{X^1\times_k X^2} \times \Delta \colon X^1 \times_k X^2 \times_k X^3 \to X^1 \times_k X^2 \times_k Z$ is trivial, as it is isomorphic to $\sL \otimes \sL^{-1} \simeq \sO$. 
This triviality implies that the diagonal $\Delta(X^3)$ is contained in $E$.
Note that  $Z$ is geometrically connected, the nonempty clopen $E$ should be $E = Z$.
Consequently, $\sN$ is trivial over $X^1 \times_k X^2 \times_k Z$. 
This precisely translates to $\sL$ satisfying the theorem of the cube.
As $X^3$ permits a section over $k$, the required $\operatorname{sq}(\sL,X^1,X^2)$ holds.
This completes the proof. \qedhere
\epf

\begin{variant}\label{variant} 
For $k$-schemes $(X^i)_{i=1}^3$ locally of finite type, geometrically (normal and integral), and an invertible sheaf $\sL$ on the cube $C$, if $X^3$ is separated, then there is a purely inseparable field extension $K/k$ such that $\mathrm{sq}(\sL_K, X^1_K,X^2_K)$ holds.
In particular, if $k$ is perfect, then $\mathrm{sq}(\sL_K, X^1_K,X^2_K)$ holds.
\end{variant}
\bpf 
By a limit argument and fpqc descent, we may assume that $k$ is perfect. 
Since the proof above essentially used the property that $X^3$ is separated and of finite type, it remains to reduce to the case when $X^3$ is of finite type.
We write $X$ as a union of increasing open subsets $X_i$, which are of finite type. 
The assertion concerns the triviality of $s(\sL)$. 
Hence, we consider the following exact sequence
\[
   0\to R^1\varprojlim H^0(X_i,\bG_m)\to H^1(X,\bG_m)\to \varprojlim H^1(X_i,\bG_m)
\]
It suffices to show that $R^1\varprojlim H^0(X_i,\bG_m)=0$, or, by \SP{0CQA}, that the system $(\GG(X_i,\sO_{X_i})^{\times})_i$ is Mittag--Leffler. 
In fact, for $i$ large enough, the injections $\GG(X_j)^{\times}/k^{\mathrm{alg}\times}\hra \GG(X_i)^{\times}/k^{\mathrm{alg}\times}$ are bijections for all $j\geq i$: these groups are of finite type due to \cite{Kah06}*{lemme~1} hence their ranks are equal when $i$ is large; by the normality of $X_j$, every $f\in \GG(X_i)^{\times}$ such that $f^n\in \GG(X_j)^{\times}$ is also invertible on $X_j$.
\epf 

\brems
\item The proof above circumvents the representability of the Picard functor, which is indispensable when $X^3$ is merely locally of finite type. 
When $X^3/k$ is proper and $\Gamma(X^3) = k$, however, an alternative perspective due to Raynaud illuminates the reason behind the theorem of the square and the necessity of taking a tensor power in \ref{power-square}.
Under these properness assumptions, $\underline{\mathrm{Pic}}_{X^3/k}$ is representable by a group scheme. The invertible sheaf $\sL$ on $X^1\times_k X^2 \times_k X^3$ induces a $k$-morphism
\[ v \colon X^1 \times_k X^2 \longrightarrow \underline{\mathrm{Pic}}_{X^3/k}. \]
Hence, the condition that $\mathrm{sq}(\sL, X^1, X^2)$ holds translates to the vanishing of the cross-difference
\[
X^1 \times_k X^1 \times_k X^2 \times_k X^2 \longrightarrow \underline{\mathrm{Pic}}_{X^3/k}
\]
given by $(x^1, y^1, x^2, y^2) \longmapsto v(x^1, x^2) - v(y^1, x^2) - v(x^1, y^2) + v(y^1, y^2)$.
Assuming $X^1$ and $X^2$ possess rational points, this vanishing is equivalent to saying $v$ is ``\emph{decomposed}'', i.e., there exist $k$-morphisms $v^i \colon X^i \to \underline{\mathrm{Pic}}_{X^3/k}$ such that $v(x^1, x^2) = v^1(x^1) + v^2(x^2)$, which corresponds to the theorem of the square.

When $X^1 \times_k X^2$ is geometrically integral, rigidity lemmata force such cross-differences to vanish if the target is an abelian variety. 
The obstruction lies entirely in the unipotent part of $\underline{\mathrm{Pic}}^0_{X^3/k}$. 
In characteristic $p>0$, any unipotent element is killed by a power $n=p^r$. 
Thus, a suitable multiple $v \mapsto nv$ eliminates the obstruction, which means that $\mathrm{sq}(\sL^{\otimes n}, X^1, X^2)$ holds.
\item The purely inseparable extension in \Cref{variant} is unavoidable.
Assume that $k$ is imperfect and $\operatorname{char}(k)=p>2$, and choose
$b\in k\backslash k^p$.
Let $G\subset \bG_{a,k}^2$ be a wound $k$-form of $\bG_{a,k}$ defined by $u-u^p=bv^p-v^{p^2}$ and let $H\subset \bG_{a,k}^2$ be the Artin--Schreier pullback of $bt^p\colon \bG_{a,k}\to \bG_{a,k}$, defined by $w+w^p+bt^p=0$.
By construction, the completion $\ov{H}\subset \mathbf{P}^2_k$ at $\infty$ is normal but over $\ov{k}$ it has a cusp so that $\Pic^0_{\ov{H}/k}$ has nontrivial unipotent part.
As $\ov{H}$ has positive genus, the Abel--Jacobi map
\[
\theta\colon H\to \Pic^0_{\ov{H}/k},\qquad x\mapsto [x-\infty]
\]
is nonconstant.
There is natural Frobenius-compatible alternating biadditive surjective pairing
\[
h\colon G\times G\to H,\qquad
((u_1,v_1),(u_2,v_2))\mapsto
\bigl(u_1u_2^p-u_1^pu_2,\; v_2u_1-v_1u_2\bigr).
\]
Then $v\ce \theta\circ h$ is not decomposed.
Let $\ov{\sL}\in \Pic(G\times G\times \ov{H})$ represent $v$,
and write $\sL\ce \ov{\sL}|_{G\times G\times H}$.
Since $v$ is not decomposed, Raynaud's criterion
(\emph{decomposed} $\Leftrightarrow$ \emph{square}) shows that
$\mathrm{sq}(\ov{\sL},G,G)$ fails.
Moreover, $s(\ov{\sL})_\xi$ is algebraically equivalent to $0$ on
$\overline X^3_\xi$, while $\ov{H}_\xi\setminus H_\xi=\{\infty_\xi\}$;
hence its nontriviality implies $s(\sL)_\xi\neq 0$.
Therefore $\mathrm{sq}(\sL, G, G)$ fails, and so does $\mathrm{cube}(\sL)$.
\erems

\bpp[Neutral components and homogeneous spaces]\label{hom-spa}
For a scheme $S$ and an $S$-group scheme $G$, recall \cite{SGA3Inew}*{exposé~IV, proposition~et~définition~6.7.3} that the $G$-homogeneous space $X$ is isomorphic to an fppf quotient of $G$ by a subgroup scheme $H\subset G$.
If the neutral component $G^{\circ}\subset G$ is representable, then $X\simeq G/H$ is equipped with an action by the group scheme $G^{\circ}$, inducing a bijection
\[
\{\text{$X/\sim$, where  $x\sim y$ if $x=g\cdot y$ for a $g\in G^{\circ}$}\}\leftrightarrow \{\text{irreducible components of $X$}\}.
\]
Hence,  every open $U\subset X$ satisfies $G^{\circ}\cdot U=X$ if and only if $U$ contains all the maximal points of $X$.
\epp

\bpp[Picard-admissible invertible sheaves]\label{gen-picard}
For a scheme $S$, an $S$-scheme $X$, and an $S$-group scheme $G$ acting on $X$, there is a $G$-action on $\Pic(X)$: for an $S$-scheme $T$ and an invertible sheaf $\sL$ on $X_T$,
\[
\text{$G(T)\times \Pic(X_T)\ra \Pic(X_T)$, $(g,\sL)\mapsto {}^{g}\sL\ce t_g^{-1}(\sL)$, where $t_g(x)=g^{-1}\cdot x$.}
\]
By construction, this action factors through $\Pic(T)$, inducing a $G$-action on the relative Picard functor 
\[
\text{$G\times \Pic_{X/S}\ra \Pic_{X/S}$,\q where\q  $\Pic_{X/S}\colon \mathbf{Sch}_{/S}\ra \Set$, \q $T\mapsto \Pic(X_T)/\Pic(T)$.}
\]
Now fix an invertible sheaf $\sL$ on $X$, then the action map above induces a natural transformation
\[
\text{$\delta_{G,\sL}\colon G\ra \Pic_{X/S}$,\qq  $g\mapsto {}^g \sL\otimes \sL^{-1}$.}
\]
Further, if $\delta_{G,\sL}$ is a homomorphism of group functors, then we say that $\sL$ \emph{is Picard-admissible}, which is equivalent to a theorem of square, see \Cref{C-equiv} below with a sketched proof.
\epp
\blem[\cite{Ray70}*{chapitre~IV, proposition~3.1}]\label{C-equiv}
For a scheme $S$, an $S$-scheme $X$, an $S$-group scheme $G$ acting on $X$, and an invertible sheaf $\sL$ on $X$, the following two assertions are equivalent
\benumr
\item\label{homo} $\sL$ is Picard-admissible, namely, $\delta_{G,\sL}\colon G\ra \Pic_{X/S}$ is a homomorphism of group functors;
\item\label{sq-holds} for the pullback $\sM$ of $\sL$ under the morphism $G\times_S G\times_S X\ra X$ by sending $(g_1, g_2,x)$ to $g_1g_2^{-1}x$,
\[
\text{ $\mathrm{sq}(\sM, G,G)$ holds. }
\]
\eenum
\elem
\bpf
The key is that,  by construction, we have $\mathrm{s}(\sM, G,G)={}^{g_1g_3}\sL\otimes {}^{g_1g_4}\sL^{-1}\otimes {}^{g_2g_3}\sL^{-1}\otimes {}^{g_2g_4}\sL$, where $(g_1,g_2,g_3,g_4)$ is the coordinate of the fourth power of $G$. 
To show that $\delta_G$ is a homomorphism, it suffices to check that ${}^{gh}\sL\otimes \sL^{-1}\simeq {}^g\sL\otimes \sL^{-1}\otimes {}^h\sL\otimes \sL^{-1}$, or ${}^{gh}\sL\otimes {}^{g^{-1}}\sL\otimes {}^{h^{-1}}\sL\otimes \sL$ is trivial, which is implied by the triviality of $s(\sM,G,G)$ by letting $g_2=g_4=1$ and $g_1=g$ and $g_3=h$.
Indeed, $s(\sM, G,G)$ is trivial in $\Pic_{X/S}$ if it is a pullback from $G\times_S G\times_S G\times_S G$, since $\Pic_{X/S}$ ignores all pullbacks from the parameter schemes functorially.
Therefore, we have \ref{sq-holds} $\Rightarrow$ \ref{homo}.
If $\delta_G$ is a homomorphism, then $s(\sM,G,G)\simeq {}^{g_1g_3g_4^{-1}g_1^{-1}g_3^{-1}g_2^{-1}g_2g_4}\sL=\sL^{\otimes 0}$ is trivial, so is a pullback from the base $G\times_S G\times_S G\times_S G$. 
\epf
\bcor\label{power3}
If $\sL$ is Picard-admissible, then for any $g\in G$, we have ${}^{g}\sL\otimes {}^{g^{-1}}\sL\simeq \sL^{\otimes 2}$ locally on $S$.
For every $s\in S$ and pair $(g,h)\in G(S)\times G(S)$, there is an open neighborhood of $s$ over which we have 
\[
\tst {}^g\sL\otimes {}^h\sL\otimes {}^{(gh)^{-1}}\sL\isoto \sL^{\otimes 3}.
\]
\ecor
% 下面是Ray70定理IV3.3.
\bprop\label{C-criterion} 
For a scheme $S$, an $S$-scheme $X$, an $S$-group scheme $G$ acting on $X$, and an invertible sheaf $\sL$ on $X$, in each of the following two cases, $\sL$ is Picard-admissible:
\benumr
\item\label{diag} $G$ is an abelian scheme, $X(S)\neq \emptyset$, and $X\times_S X$ is its only clopen containing $\Delta_{X/S}(X)$.
\item\label{Lt-C} $S$ is locally Noetherian, $X$ has an $S$-section, $G$ is flat finite type over $S$ such that $\sO(G_s)=k(s)$ for all $s\in S$, and for every $t\in S$ as the image of a maximal point of $X\times_S X$, $\sL_{t}$ is Picard-admissible.
\eenum
Suppose that $S$ is normal integral, $G$ and $X$ are of strict (N)-type (resp, (N)-type) over $S$ with geometrically integral generic fibres, and for each $s\in S$ (resp, such that $\mathrm{wdim}\,\sO_{S,s}=1$), either
\benuma
\item\label{Pic-comp-Gs-int-X-section}  $G_s$ is integral and $X$ has an $S$-section, or
\item\label{Pic-comp-Xs-geomint}  $X_s$ is geometrically integral.
\eenum
There is an integer $n>0$ such that $\sL^{\otimes n}$ is Picard-admissible.

% \item\label{Prufer-C} $S$ is integral and normal, $G$ and $X$ are of (N)-type over $S$, such that either
% \benum
% \item $X$ has an $S$-section and $G_s$ is integral for every $s\in S$,
% \item $X_s$ is geometrically integral for every $s\in S$.
% \eenum
% \eenum
\eprop
\bpf
The cases \ref{diag} and \ref{Lt-C} are \cite{Ray70}*{chapitre~IV, théorème~3.3 1) and 2)}.  
We prove the next assertions.
Recall the Picard-admissibility, we need to show that there is an integer $n>0$ such that the pullback $\sM^{\otimes n}$ of $\sL^{\otimes n}$ along the morphism $G\times_S G\times_S X\ra X$ via $(g_1,g_2,x)\mapsto g_1g_2^{-1}x$ makes $\mathrm{sq}(\sM^{\otimes n}, G,G)$ hold. 
If (2) is satisfied, then it suffices to let $X^1=X^2=G$ and $X^3=X$ in  \Cref{power-cube-square}\ref{power-square} to conclude that $\mathrm{sq}(\sM^{\otimes n},G,G)$ holds.
If (1) is satisfied, then we regard $X$ as $X^3$ at the end of \S \ref{thm-sq} by the presence of the $S$-section of $X$ to reduce $\mathrm{sq}(\sM^{\otimes n},G,G)$ to $\mathrm{cube}(\sM^{\otimes n})$.
Now, we change the order of $G,G$, and $X$ by letting $X^1=X$ and $X^2=X^3=G$.
By \cite{EGAIV2}*{corollaire~4.5.14}, the scheme $X^3=G$ satisfies the conditions of \Cref{power-cube-square}.
Also, $G$ and $X$ are of (strict) (N)-type  over $S$ with geometrically integral generic fibres, so \Cref{power-cube-square}\ref{power-square} implies that the desired $\mathrm{cube}(\sM^{\otimes n})$ holds.
\epf

\section{Criteria for ampleness}\label{sec-crit}
This section translates the Picard-admissibility established earlier into criteria for ampleness. Bypassing Noetherian hypotheses, we explicitly construct ample sheaves directly from group orbit boundaries.
\bpp[Projectivity and group actions]\label{proj-gp}
Let $f\colon X\ra S$ be a quasi-compact and quasi-separated morphism of schemes. 
For an invertible $\sO_X$-module $\sL$ and the quasi-coherent graded $\sO_S$-module $\cB_{\bullet}\ce f_{\ast}(\bigoplus_{n\geq 0}\sL^{\otimes n})$,  the map $\id \colon \cB_{\bullet}\ra \cB_{\bullet}$ induces a canonical homomorphism of graded $\sO_X$-algebras
\[
 \tst  \psi\colon f^{\ast}(\cB_{\bullet})\ra \bigoplus_{n\geq 0}\sL^{\otimes n}.
\]
Note that $\Proj (\bigoplus_{n\geq 0}\sL^{\otimes n})\cong \Proj(\sO_X[t])=X$, the homomorphism $\psi$ induces the following $S$-morphism 
\[
r_{\sL}\colon X_{\sL}\ra P\ce \Proj(\cB_{\bullet}),
\]
where $X_{\sL}\subset X$ is  locally the complementary open of $V_+(\psi(f^{\ast}(\cB_{\bullet})_+))$, see \cite{EGAII}*{2.8.1}.
In particular, $X_{\sL}$ is the largest open over which $\psi$ is surjective on large degrees. 
Geometrically speaking, if $X_{\sL,i}\subset X$ is the open locus where $\sL^{\otimes i}$ is $S$-relatively generated by global sections for each $i\in \mathbb{Z}_{>0}$, then we have
\[
  \tst  \x{$X_{\sL}=\bigcup_{i>0}X_{\sL,i}$,\q where $X_{\sL,i}=\bigcup_{s\in \GG(X,\sL^{\otimes i})}X_{s}$}
\]
and $X_{\sL,i}$ is the largest open of $X$ over which $f^{\ast}f_{\ast}\sL^{\otimes i}\surjects \sL^{\otimes i}$ is surjective.
In particular, by \cite{EGAIII1}*{proposition~1.4.15} and \cite{EGAIV1}*{1.7.21}, the formation of $X_{\sL,i}$ commutes with flat base change $S\pr\ra S$.

\textbf{Setup.} Let $G$ be an $S$-group scheme acting on $X$.
We assume that $G$ is universally open over $S$, so that $G$-orbit of every open  $U\subset X$ is open in $X$.
Let $\sL$ be Picard-admissible with the $G$-action.
\epp
%\blem
%For a scheme $S$, an $S$-scheme $X$, an irreducible $S$-group scheme $G$ acting on $X$, a point $x\in X(S)$, a element $g\in G(S)$, and an open subscheme $U\subset X$ such that one of the orbit $G\cdot x$ factors through $U$, there exists a nonempty open subscheme $W\subset G$ such that 
%
%\elem

\bprop\label{big-prop}
Under the Setup, suppose that $f$ is of finite type, the $S$-group $G$ is flat finitely presented with connected $S$-fibres, and $X=\bigcup_{i>0} X_{\sL,i}$.
Then, for every $g\in G(S)$ there exist
\[
\text{open subschemes $U_g$ and $U_g\pr$ of $P$ containing $r_{\sL}(X)$\q  and\q  an $S$-isomorphism $\tau_g\colon U_g\isoto U_g\pr$}
\]
such that the following diagram is commutative, where $t_g$ is the translation by $g$:
\[
\begin{tikzcd}
X \arrow[r] \arrow[d, "t_g"'] \arrow[rr, "r_{\sL}", dashed, bend left] & U_g \arrow[r] \arrow[d, "\tau_g"'] & P\, \\
X \arrow[r] \arrow[rr, "r_{\sL}"', dashed, bend right]                 & U_g^{\prime} \arrow[r]             & P\,.
\end{tikzcd}
\]
In particular, the following assertions hold
\benumr
\item The set $X^{\mathrm{fl},\sL}\ce\{x\in X\,|\, \text{$r_{\sL}$ is flat at $x$}\}$ satisfies $G\cdot X^{\mathrm{fl},\sL}=X^{\mathrm{fl,\sL}}$;
\item The set $X^{\mathrm{qf},\sL}\ce \{x\in X\,|\, \text{$r_{\sL}$ is quasi-finite at $x$}\}$ is open in $X$ and $G\cdot X^{\mathrm{qf},\sL}=X^{\mathrm{qf},\sL}$;
\item\label{g-stab} The largest open $P_0\subset P$ such that $r_{\sL}|_{P_0}$ is an open immersion has $G$-saturated preimage 
\[
G\cdot r_{\sL}^{-1}(P_0)=r_{\sL}^{-1}(P_0)
\]
 and restriction of $\sL$ on every quasi-compact open subscheme of $ r_{\sL}^{-1}(P_0)$ is $S$-ample.
Further, if $X$ is $S$-separated, then we have $ r_{\sL}^{-1}(P_0)=X^{\mathrm{qf},\sL}$.
\eenum
\eprop
\bpf
See \cite{Ray70}*{chapitre~V, proposition~3.1 and corollaire~3.6}.
\epf

\blem\label{Ray70-V2.4} 
For a scheme $S$, an $S$-flat group scheme $G$ with connected fibres acting on a quasi-compact quasi-separated $S$-scheme $X$, and a Picard-admissible invertible sheaf $\sL$ on $X$, 
\benumr
\item\label{contain} we have $G\cdot X_{\sL,i}\subset X_{\sL,2i}$ for all $i\in \mathbb{Z}_{>0}$;
\item\label{open-stable} the open $X_{\sL}\subset X$ is stable under $G$-action;
\item\label{entire-intersect} $X_{\sL}=X$ if and only if for every $x\in X$, the intersection $X_{\sL}\cap G\cdot x$ is nonempty;
\item\label{generated-globally} if $X$ is quasi-compact, then \ref{entire-intersect} holds if and only if $\sL$ is $S$-semiample (that is, some power $\sL^{\otimes n}$ is generated by its global sections relative to $S$).  
\eenum
\elem
\bpf
The assertions \ref{open-stable}--\ref{generated-globally} follow from \ref{contain} directly.
The proof of \ref{contain} is given in \cite{Ray70}*{chapitre~V, proposition~2.1} and we sketch it briefly.
For every $g\in G(S)$, the translation ${}^{g}\sL^{\otimes i}$ is generated by global sections on $g\cdot X_{\sL,i}$.
If $g\cdot X_{\sL,i}\cap g^{-1}\cdot X_{\sL,i}\neq \emptyset$, then by  ${}^{g}\sL^{\otimes i}\otimes {}^{g^{-1}}\sL^{\otimes i}\simeq \sL^{\otimes 2i}$ (\Cref{power3}), the power $\sL^{\otimes 2i}$ is generated by its sections on $g\cdot X_{\sL,i}\cap g^{-1}\cdot X_{\sL,i}$ locally on $S$, to the effect that
\[
g\cdot X_{\sL,i}\cap g^{-1}\cdot X_{\sL,i}\subset X_{\sL,2i}.
\]
By this observation, it suffices to show that each $x\in G\cdot X_{\sL,i}$ lies in some $g\cdot X_{\sL,i}\cap g^{-1}\cdot X_{\sL,i}$ for a $g\in G(S)$. 
By taking a flat base change from a local ring of $G$, we may assume that $S$ is local, so we are reduced to finding such $g$ in the closed fibre $G_s$.
Hence everything is over a field, then the connectedness of $G_s$ (so it is geometrically irreducible) implies the existence of such $g$.
\epf
%V3.9
\bprop\label{ample-criterion}
Let $R$ be a normal ring and $G$ a flat finite type $R$-group scheme with connected fibres.
For a topologically Noetherian $R$-scheme $X$ of strict (N)-type (resp., (N)-type) and an invertible $\sO_X$-module $\sL$, if $G$ is $S$-smooth (resp., $G_s$ is smooth if $\mathrm{wdim}\, \sO_{S,s}\leq 1$), then we have equivalences:
\benum
\item\label{semiample-cr} $\sL$ is semiample iff every orbit $G\cdot x$ meets some $X_f$ for $f\in \GG(X,\sL^{\otimes n})$;
\item\label{ample-cr} $\sL$ is ample iff every orbit $G\cdot x$ meets a quasi-affine $X_f$ for $f\in \GG(X,\sL^{\otimes n})$.
\eenum
\eprop
\bpf 
Since the implications ``$\Rightarrow$'' are clear for (a) and (b), it remains the directions ``$\Leftarrow$''.
\Cref{normality}\ref{N-normal} and \ref{weak-N-normal} yield the normality of $X$ respectively.
Every quasi-compact open of $X$ has finitely many irreducible components, so by \SP{0357}, we may replace $X$ by each of its normal integral connected components (with generic point $\xi$).
Denote $S\ce \Spec R$ and $\eta$ the generic point.

We first consider the preliminary case when $X_{\eta}$ is geometrically integral and $X$ has an $S$-section. 
As \cite{EGAIV2}*{corollaire~4.5.14} implies that $G$ and $X$ have geometrically integral generic fibres, the conditions in \Cref{C-criterion} (1) are satisfied, so there exists an integer $m>0$ such that $\sL^{\otimes m}$ is Picard-admissible. 
Therefore, we may replace $n$ by $mn$ to assume that $\sL$ is Picard-admissible. 
For the case (a), we apply \SP{047M} and \Cref{Ray70-V2.4}\ref{generated-globally} to conclude. 
For the case (b), by \SP{01Q3}(1)$\Leftrightarrow$(2), it suffices to show that for all homogeneous elements $b\in \GG(X,\sL^{\otimes \bullet})_{+}$, the open subsets $X_b$ cover $X$ and the canonical morphism $r_{\sL}\colon X\ra P$ is an open immersion.
For the largest open subset $P_0\subset P$ over which $r_{\sL}$ is an open immersion, the $G$-stability of $r^{-1}_{\sL}(P_0)$ (\Cref{big-prop}\ref{g-stab}) and the assumption reduce us to showing that every $X_b$ is contained in $r^{-1}_{\sL}(P_0)$.
Equivalently, it suffices to show that $r_{\sL}|_{X_b}$ is an open immersion.
Recall \cite{EGAII}*{proposition~3.7.3} that $r_{\sL}^{-1}(D_+(b))=X_{b}$ and $\GG(X_b,\sO_X)\simeq \GG(D_+(b),\sO_P)$, where $D_+(b)\subset P$ is the affine open determined by $b$.
Hence, the morphism $r_{\sL}|_{X_b}\colon X_b\ra D_+(b)=\Spec \GG(X_b, \sO_X)$ is an open immersion because $X_b$ is quasi-affine. 

Now, we establish the general case by using above resolved case iteratively.
Recall the notation in \Cref{reduce-geonum-irrcomp}.
Our induction hypothesis is that the assertion  holds for all $X$ with $\rho(X)\leq N$ for a fixed integer $N\geq 0$.
Thus, the goal is to show that under this hypothesis, the assertion holds when $\rho(X)=N+1$.
Note that the case when $N=0$ is trivial because $X=\emptyset$.
Now assume that $N\geq 1$.

Replacing $S$ by its open subscheme, the image of $X$, we may assume that $X\ra S$ is fpqc.
If $X\times_S X$ is not connected, then we exploit \Cref{reduce-geonum-irrcomp} so that $\rho(X_i)\leq \rho(X)-1=N$, hence the induction hypothesis applies to each $X_i$.
In particular, the assertion holds for $X\times_S X$.
On the other hand, if $X\times_S X$ is connected, then the analysis in the proof of the claim already shows that $X\times_S X$ has geometrically integral generic fibre over $X$, so the preliminary case applies to $X\times_S X$.
By fpqc descent, the assertion for $X\times_S X$ descends to $X$, hence the assertion holds when $\rho(X)=N+1$, as desired.
\epf
\bcor\label{V3.8}
For an affine scheme $S$, an $S$-quasi-separated finite type scheme $X$, an $S$-flat finitely presented group scheme $G$ with connected fibres acting on $X$,  a Picard-admissible invertible $\sO_X$-module $\sL$ is $S$-ample if and only if $X$ is covered by the $G$-orbit of $U=\bigcup_{s\in \GG(X,\sL^{\otimes n})} X_s$ where $X_s$ are quasi-affine.
\ecor
\bpf
If $\sL$ is $S$-ample, then by \cite{EGAII}*{théorème~4.5.2} every $x\in X$ has an affine open neighborhood of the form $X_s$ for some $s\in \GG(X,\sL^{\otimes n})$. 
In particular, we have $U=X$, hence ``$\Rightarrow$'' is proved.

Now, we show ``$\Leftarrow$''. 
First, using \Cref{Ray70-V2.4}\ref{entire-intersect}, we have $X=X_{\sL}$.
By \Cref{big-prop}\ref{g-stab}, the preimage $r_{\sL}^{-1}(P_0)$ is $G$-stable.
Furthermore, to prove that $\sL$ is $S$-ample, it suffices to show that $r_{\sL}^{-1}(P_0)=X$. Indeed, if so, then $X$ is covered by quasi-compact open subschemes on which $\sL$ is ample, which implies $\sL$ is ample on $X$ since $X$ is quasi-compact.
Thus we are reduced to showing that $U\subset r_{\sL}^{-1}(P_0)$, or, for every $\sigma\in \GG(X,\sL^{\otimes n})$ such that $X_{\sigma}$ is quasi-affine, the restriction $r_{\sL}|_{X_{\sigma}}\colon X_{\sigma}\hra P$ is an open immersion.

Denote $P_{\sigma}\ce D_{+}(\sigma)=\Spec ((\cB_{\bullet})_{\sigma})_0$, then $r_{\sL}$ restricts to a morphism $ X_{\sigma} \to P_{\sigma}$.
This induces an isomorphism $\GG(P_{\sigma},\sO_{P})\xrightarrow{\sim} \GG(X_{\sigma}, \sO_X)$ by sending $t/\sigma^k$ where $s\in \GG(X,\sL^{\otimes nk})$ to the unique function $f\in \GG(X_{\sigma},\sO_X)$ such that $t|_{X_{\sigma}}=f\cdot \sigma^k|_{X_{\sigma}}$. 
Since $P_{\sigma}$ is affine, we can identify $P_{\sigma}$ with $\Spec(\GG(X_{\sigma}, \sO_X))$.
By \cite{EGAII}*{proposition~5.1.2}, the canonical morphism from a quasi-affine scheme to the spectrum of its global sections is an open immersion.
Therefore, the morphism $r_{\sL}|_{X_{\sigma}}$ is an open immersion.
\epf

\bcor[\cite{Ray70}*{chaptre~VI, théorème~2.1}]\label{homo-ample-crit}
Let $S$ be an affine scheme and $G$ an $S$-flat locally finitely presented group scheme.
Let $Y$ be a $G$-homogeneous space locally of finite type over $S$.
Suppose that the neutral component $G^{\circ}\subset G$ is representable and there is a $G^{\circ}$-stable quasi-compact quasi-separated open $X\subset Y$.
For a Picard-admissible invertible $\sO_X$-module $\sL$ with respect to the $G^{\circ}$-action, then   
\benumr
\item\label{semiample-crit} $\sL$ is semiample iff any maximal point of fibres of $X$ lies in an $X_{\sigma}$ for $\sigma\in \GG(X,\sL^{\otimes n})$;
\item\label{ample-crit} $\sL$ is ample iff any maximal point of fibres of $X$ lies in a quasi-affine $X_{\sigma}$ for $\sigma\in \GG(X,\sL^{\otimes n})$.
\eenum
\ecor
\bpf
For the case \ref{semiample-crit}, by the quasi-compactness of $X$,  the sheaf $\sL$ is semiample if and only if $X=X_{\sL,n}$ for some $n>0$, if and only if $X=X_{\sL}$, and by the $G^{\circ}$-stability of $X_{\sL}$, if and only if $G^{\circ}\cdot X_{\sL}=X$, hence by the end of \S\ref{hom-spa}, is equivalent to the existence of a desirable $X_{\sigma}$.
For the case \ref{ample-crit}, we take $U$ as the union of all quasi-affine open $X_{\sigma}\subset X$ for all $\sigma\in \GG(X,\sL^{\otimes n})$ for some positive integer $n$, and then apply the criterion for ampleness \Cref{V3.8} to acquire the equivalence between the ampleness of $\sL$ and the assertion $X=G^{\circ}\cdot U$, which is equivalent to the existence of a desirable $X_{\sigma}$ by the end of \S\ref{hom-spa}.
\epf
\blem\label{separated} 
For a scheme $S$, let $G$ be an $S$-group scheme that is universally open, locally of finite type, with connected $S$-fibres and let $X$ be  a  locally finite type $S$-scheme equipped with a $G$-action.
Then 
\[
   \x{$X$ is $S$-separated\q iff \q there is an $S$-separated open $U\subset X$ such that $G\cdot U=X$.}
\]
\elem
\bpf 
We use valuative criterion for separatedness. 
Since all the conditions are preserved under $S$-base change, it suffices to assume that $S$ is the spectrum of a valuation ring $V$ with fraction field $K$.
Since the $S$-separatedness of $X$ is insensitive to its nilpotent structure, we may assume that $X$ and $G$ are reduced. 
Then by \Cref{flat-equiv}, $G$ is $S$-flat. 
The algebraic closure $\ov{K}$ contains $V$ so by \SP{00IA} we may assume that ${K}$ is algebraically closed, hence so is the residue field $\kappa\ce V/\fm_V$. 
By \Cref{raygru}, $G$ is locally of finitely presentation over $V$.
Now, we consider two $V$-sections $s_1$ and $s_2$ of $X$ such that $s_1|_K=s_2|_K$ and the goal is to prove that $s_1=s_2$.
If $s_1$ and $s_2$ are contained in $U$, then we are done.
Hence, it suffices to exploit the $G$-action to transport $s_i$ inside $U$.
First, consider two maps $\phi_i\colon G_{\kappa}\ra X_{\kappa}$ by sending $g$ to $g\cdot (s_i)_{\kappa}$.
Since $G$-translations of $U$ cover $X$, the preimages $\phi_i^{-1}(U)$ are nonempty.
Besides, the Cohen-Macaulay locus $\mathrm{CM}(G_{\kappa})$ is open in $G_{\kappa}$ \SP{00RF}.
Note that the connectedness of $G_{\kappa}$ implies the irreducibility of $G_{\kappa}$ \cite{SGA3Inew}*{exposé~VI, corollaire~2.4.1}, so we have 
\[
   W\ce \phi_1^{-1}(U)\cap \phi_2^{-1}(U)\cap \mathrm{CM}(G_{\kappa})\neq \emptyset. 
\]
As $\kappa$ is algebraically closed, there is an element $g_0\in W$.
Thus, by \cite{EGAIV4}*{proposition~17.16.1}, there exists a closed subscheme $T\subset G$ that is fppf finite over $V$ and passing through $g_0$.
Taking a valuation ring $V\pr$ dominating the local ring of $T$ at $g_0$, we obtain a flat local homomorphism $V\hra V\pr$, so it is an fpqc cover.
Base changing everything to $V\pr$, we extend $g_0$ to a $V\pr$-section $g$ of $G_{V\pr}$.
As $U_{V\pr}$ is separated, two sections are equal $g\cdot (s_1)_{V\pr}=g\cdot (s_2)_{V\pr}$, hence  $(s_1)_{V\pr}$ and $(s_2)_{V\pr}$ coincide.
Every algebraic space over a scheme is an fpqc sheaf by Gabber \SP{0APL}, we conclude that $s_1=s_2$, so $X$ is $S$-separated.
\epf
% 以下是V3.10
\bthm\label{constr-ample}
Let $S$ be a locally (coherent, topologically Noetherian) normal scheme and $G$ an $S$-flat group of finite type with connected fibres.
For an $S$-smooth scheme $X$ with $G$-action and an $S$-quasi-affine open $U$, let $(D_i)_{i\in I}$ be one-codimensional irreducible components of $G\cdot U - U$. 
If $G_s$ is smooth whenever $\mathrm{wdim}\,\sO_{S,s}\leq 1$, then any cycle $D\ce \sum_{n_i\in \bZ_{>0}}n_i D_i$  is an effective divisor, and $\sL\ce \sO_X(D)$ is $S$-ample.
\ethm
\bpf
Since $U$ is $S$-quasi-affine, by \SP{01SL} it is quasi-compact and separated and the orbit $G\cdot U$ is quasi-compact.
We may assume that $X=G\cdot U$, which is $S$-separated due to \Cref{separated}.
The problem is local on $S$, so we may assume that $S$ is affine to prove the ampleness of $\sL$.
Note that each $S$-fibre of $G$ is geometrically irreducible, so $U$ contains all the maximal points of $S$-fibres of $X$.

\textbf{Step 0.} If $I=\emptyset$, then $\sL=\sO_X$.  
The Hartogs's extension \cite{GL24}*{Theorem~2.20} implies that $j_{\ast}\sO_U\simeq \sO_X$, where $j\colon U\hra X$ is the open immersion, and similar for every  $G$-translation of $U$.
Every $x\in X$ is contained in a $G$-translation of $U$, whose $S$-quasi-affineness yields an affine open $U_b$ where $b\in \GG(U, \sO_U)$.
Then $b$ extends to $\wt{b}\in \GG(X,\sO_X)$ such that $x\in X_{\wt{b}}\subset X$.
The $S$-separatedness of $X$ and the affineness of $U_b$ imply that $U_b\hra X_{\wt{b}}$ is affine.
However, we have $\GG(X_{\wt{b}},\sO_X)\simeq \GG(U_b, \sO_X)$ so $U_b=X_{\wt{b}}$.
As a result, $X_{\wt{b}}$ is an affine open containing $x$, so $\sL$ is ample (\cite{EGAII}*{définition~4.5.3}).

\textbf{Step 1.} If $I\neq \emptyset$, then the cycle $D$ is a relative effective  divisor on $X$.
To see this, we factorize $j$ as $U\hra U^{\mathrm{aff}}\hra X$, where $U^{\mathrm{aff}}=\un{\Spec}_{X}j_{\ast}\sO_U$ is the affinization such that $\depth \sO_{X,z}\geq 2$ for every $z\in U^{\mathrm{aff}}\backslash U$. 
Then the open immersion $U^{\mathrm{aff}}\hra X$ is affine and $D\subset X\backslash U^{\mathrm{aff}}$, so \Cref{Ramanujam--Samuel} applies.

\textbf{Step 2.} To verify that $\sL=\sO_X(D)$ is $S$-ample, we leverage the criterion for ampleness \Cref{ample-criterion} \ref{ample-cr} to reduce us to showing that for every $x\in X$, there is $f\in \GG(X,\sL^{\otimes n})$ such that $X_f$ is quasi-affine and $X_f\cap G\cdot x\neq \emptyset$.  
Since $X=G\cdot U$ and $U$ is quasi-affine, we may assume that $x\in U$.
The faithfully flat descent for ampleness and the d\'evissage in the proof of \Cref{ample-criterion}
reduces us to the case where $X_{\eta}$ is geometrically integral and $X$ has an $S$-section.
Therefore, replacing $\sL$ by a larger power as in \Cref{C-criterion}\ref{Pic-comp-Xs-geomint}, we may assume that $\sL$ is Picard-admissible. 

\textbf{Step 3.} Recall \S\ref{proj-gp} that there exist an open subscheme $X_{\sL}\subset X$ and an $S$-morphism
\[
  \text{$r_{\sL}\colon X_{\sL}\ra P$.} 
\]
Since $D$ is effective, by \SP{01X0}, there exists $\sigma\in \GG(X,\sL)$ such that $X\backslash \Supp(D)=X_{\sigma}$.
By construction of $D$, we have $U\subset X_{\sigma}$ so $G\cdot X_{\sigma}=X$ and \Cref{Ray70-V2.4}\ref{entire-intersect}\ref{generated-globally} implies that $X_{\sL}=X$.
On the other hand,  we have a Hartogs's extension $\GG(X_{\sigma},\sO_X)\simeq \GG(U,\sO_X)$.
For the morphism $r_{\sL}\colon X\ra P$, by \cite{EGAII}*{proposition~3.7.3}, we have $X_{\sigma}=r_{\sL}^{-1}(D_+(\sigma))$ and $\GG(D_+(\sigma),\sO_P)\simeq \GG(X_{\sigma},\sO_X)$.
Consequently, by composition, we obtain an isomorphism $\GG(D_+(\sigma),\sO_P)\simeq \GG(U,\sO_X)$.
Note that $U$ is $S$-quasi-affine, thus \cite{EGAII}*{proposition~5.1.2} implies the desirable result that $r_{\sL}|_U\colon U\hra P$ is an open immersion.

\textbf{Step 4.} For an $x\in U$, since $P$ is a projective space, there exists $f\in \GG(X,\sL^{\otimes n})$ such that $D_+(f)\subset P$ is open affine and $r_{\sL}(x)\subset D_+(f)\subset r(U)$.
By Step 3, $r_{\sL}|_U$ is an open immersion, so the preimage $X_{f}\ce r_{\sL}^{-1}(D_+(f))$  satisfies $X_{f}\cap U\simeq D_+(f)$, which is affine and contains $x$.
The open immersion $X_{f}\cap U\hra X_{f}$ is affine thanks to \cite{EGAII}*{corollaire~1.6.3}, then $X_f\cap U=(X_f\cap U)\aff=X_f$ because $\GG(X_f,\sO_X)\simeq \GG(X_f\cap U, \sO_X)$.
Consequently, we obtain $X_{f}\subset U$ as a quasi-affine neighborhood of $x$.
\epf
By a limit argument, the smooth case of \Cref{constr-ample} leads to \cite{Ray70}*{chapitre~V, corollaire~3.14}.
\bcor
Let $R$ be a normal domain and $G$ an $R$-smooth group with connected fibres acting on an (N)-type $R$-scheme $X$.
If $X$ has an $R$-quasi-affine open $U$ such that $G\cdot U=X$, then $X$ is quasi-projective.
\ecor
\bcor\label{VI2.3} 
For a locally (coherent, topologically Noetherian) normal scheme $S$, an $S$-smooth group $G$, a $G$-homogeneous space $Y$ with opens $U\subset X\subset Y$ such that $U$ is $S$-quasi-affine and fibrewise dense in $X$, and all irreducible components $(D_i)_{i\in I}$ of $X\backslash U$ of codimension one in $X$, 
then for any $n_i\in \mathbb{Z}_{>0}$
\[
   \tst \x{$D\ce \sum_{i\in I}n_iD_i$ is an effective divisor, and $\sL\ce \sO_{X}(D)$ is $S$-ample.}
\]
\ecor
\bpf
The smoothness of $G$ by \cite{SGA3Inew}*{exposé~VI$_{\text{B}}$, théorème~3.10(iv)} yields the representability of the neutral component $G^{\circ}$, hence we may consider the open subsets $X \subset G^{\circ}\cdot U\subset Y$ (where the first inclusion is due to the end of \S\ref{hom-spa}).
Subsequently, let $(F_j)_{j\in J}$ be the irreducible components of $G^{\circ}\cdot U\backslash X$ of codimension one in $G^{\circ}\cdot U$ (and are $S$-flat in the case (ii)) and let $D\pr\ce \sum_{j\in J}F_j$, then \Cref{constr-ample} implies that $D+D\pr$ is an effective Cartier divisor on $G^{\circ}\cdot U$ such that $\sO_{G^{\circ}\cdot U}(D+D\pr)$ is $S$-ample. 
The invertible sheaf $\sO_X(D)$ is the restriction of the $S$-ample sheaf $\sO_{G^{\circ}\cdot U}(D+D\pr)$ on $X$, hence is $S$-ample.
\epf

The following \Cref{finite-contain} generalizes \SP{09NN} and
will be used to find an open subscheme containing all the maximal points of fibres of schemes flat over Pr\"ufer rings, see \Cref{cor-prufer-qproj}.
\bprop\label{finite-contain}
Let $X$ be a separated scheme whose every quasi-compact open subset has a finite number of irreducible components.
For a finite set $S\ce \{x_1,\cdots, x_r\}$ such that every $x_i$  is dominated by only finitely many valuation rings of finite ranks, there exists an affine open of $X$ containing  $S$.
\eprop
\bpf
We may assume that $X$ is integral with function field $K$. Let $X^{\nu}$ be the normalization. 
Suppose that there is already an affine neighborhood $U$ of $S$. 
Then by using \SP{05YU} and the prime avoidance, we are reduced to the case when $X$ is normal.
By the separatedness, each $x_i$ is dominated by a unique valuation.
Recall \cite{Bou98}*{chapitre~VI, \S7, propositions~1, 2} that the intersection $R\ce \bigcap_{i=1}^r\sO_{X,x_i}$ is a semilocal Prüfer domain. 
As $\sO_{X,x_i}$ are finite dimensional, there is an affine open covering of $T\ce \Spec R$ by principal open subsets, gluing into a morphism $T\to X$, so that we conclude by \SP{09NI}. 
\epf 

\bcor\label{cor-prufer-qproj} 
For a semilocal Pr\"ufer scheme and an $S$-smooth finite type group scheme $G$, every $S$-separated $G$-homogeneous space is $S$-quasi-projective.
\ecor
\bpf
There are finitely many maximal points of $S$-fibres of $X$, so \Cref{finite-contain} and the $S$-separatedness of $X$ yield an affine open $U\subset X$ containing all the maximal points of $S$-fibres of $X$. 
In particular, we have $G^{\circ}\cdot U=X$ due to the end of \S\ref{hom-spa}, so \Cref{VI2.3} implies that $X$ is $S$-quasi-projective.
\epf
\bcor
Let $S$ be a locally coherent normal scheme and a homogeneous space $Y$ under a smooth $S$-group $G$.
Every open $X$ of $Y$ that is surjective over $S$ with connected fibres is locally $S$-quasi-projective.
\ecor
\bpf
We may assume that $S$ is affine and take an affine open $U\subset X$ whose open image is denoted by $W\subset S$. 
It suffices to prove that the surjective morphism $X_W\ra W$ is quasi-projective.
Since $X$ has $S$-connected fibres, which are irreducible due to the $S$-smoothness of $X$, hence $U$ contains all the maximal points of $S$-fibres of $X$. 
It remains to apply \Cref{VI2.3} to finish the proof.
\epf
Finally, we acquire the local quasi-projectivities as \cite{Ray70}*{chapitre~V, 3.14 and chapitre~VI, 2.4}.
\bcor
For a normal scheme $S$ and an $S$-smooth group scheme $G$ with connected fibres, a locally finite type $S$-scheme equipped with a $G$-action, assume that one of the following conditions holds
\benumr
\item\label{qp-i}  $X$ is locally finitely presented over $S$ with a quasi-affine open $U\subset X$ such that $G\cdot U=X$;
\item\label{qp-ii}  $X$ is a $G$-homogeneous space.
\eenum
Then $X$ is locally $S$-quasi-projective. In particular, if $S$ is affine integral, then $X$ is $S$-quasi-projective.
\ecor
\bpf
For the case \ref{qp-i}, it suffices to invoke the standard limit argument to reduce to the case when $S$ is Noetherian, and then apply \Cref{constr-ample}.
The case \ref{qp-ii} follows from \ref{qp-i} by localization on $S$ with a limit argument, where we used the facts that $X$ is locally finitely presented over $S$, and that all $S$-fibres of $X$ are irreducible to find the quasi-affine $U\subset X$ in \ref{qp-i}.
\epf

\section{Extension of ample invertible sheaves}\label{sec-ext-ample}
Building on our valuative framework, this section establishes criteria for extending generic polarizations over integral Prüfer bases. We demonstrate how generic ampleness uniquely and automatically extends to global S-ampleness, circumventing the failure of standard spreading-out arguments.
\blem\label{affine-contain-ample-power-extend-ample} 
For an integral Pr\"ufer scheme $S$ with generic point $\eta$, an $S$-smooth group scheme $G$, a $G$-homogeneous space $Y$ with a $G^{\circ}$-stable $S$-quasi-compact open $X\subset Y$, and an $S$-ample invertible sheaf  $\sL_{\eta}$ on $X_{\eta}$, if $X$ has an $S$-affine open $U$  containing all the maximal points of $S$-fibres of $X$, then
\[
\text{there is $n\in \mathbb{Z}_{>0}$ such that $\sL_{\eta}^{\otimes n}$ extends to an $S$-ample invertible sheaf $\sM$ on $X$.}
\]
\elem
\bpf
Since $\sL_{\eta}$ is ample on $X_{\eta}$, by \SP{09NV}, there is $f\in \GG(X_{\eta},\sL^{\otimes n})$ such that the affine open $(X_{\eta})_{f}\subset U\cap X_{\eta}$ contains all the maximal points of $X_{\eta}$. 
Hence, the pair $(\sL_{\eta}^{\otimes n},f)$ corresponds to an effective divisor $D_{\eta}$ with $\Supp(D_{\eta})=Z(f)\ce\{f=0\}$.
Since $(X_{\eta})_{f}\subset U_{\eta}$, we get $\Supp(D_{\eta})\supset X_{\eta}\backslash U_{\eta}$.
The schematic closure of $Z(f)$ in $X$ by \Cref{Ramanujam--Samuel} defines an effective divisor $D$ on $X$.
We claim that $X\backslash U\subset \Supp(D)$.
The hypothesis on $U$ and \Cref{separated} lead to the $S$-separatedness of $X$.
As $U$ is affine over $S$, by \SP{01SG}, the open immersion $U\hra X$ is affine.
Note that $X\backslash U$ fibrewise satisfies the condition in \cite{EGAIV4}*{corollaire~21.12.7}, hence $X\backslash U$ is of pure codimension one in $X$ and is the schematic closure of its generic fibre $(X\backslash U)_{\eta}$ in $X$.
The inclusion $X_{\eta}\backslash U_{\eta}\subset \Supp(D_{\eta})$ then implies that $X\backslash U\subset \Supp(D)$ and $X\backslash \Supp(D)\subset U$.
Therefore, since $U$ is affine, $X\backslash \Supp(D)$ is quasi-affine over $S$.
The invertible $\sO_X$-module $\sM\ce \sO_{X}(D)$ is $S$-ample by \Cref{VI2.3}.
\epf
\blem\label{one-dim-semilocal} 
For a valuation ring $V$ with fraction field $K$ and a finite field extension $L/K$,  the integral closure $R$ of $V$ in $L$ is semilocal Pr\"uferian, each fibre of $\Spec R\to \Spec V$ is finite and $\dim R=\dim V$.
\elem
\bpf 
By \cite{Bou98}*{chapitre~VI, \S8, \textnumero3, remarque}, the integral closure $R$ of the valuation ring $V$ in $L$ is semilocal Pr\"uferian.
It follows from \SP{00OK} that $R$ also has Krull dimension $\dim V$.
It remains to see that each fibre of $\Spec R\to \Spec V$  is a finite set.
As normalization commutes with localizations, it suffices to consider the closed fibre.
Each points of the closed fibre corresponds to an extension of the valuation $v\colon K\to \GG$ such that the valuation ring of $v$ is $V$.
Since $L/K$ is finite, there are only finitely many extensions of $v$ to $L$. 
Therefore, each fibre of $\Spec R\to \Spec V$ is a finite set.  
\epf
Now we acquire a descent of ample invertible sheaves over valuation rings (\emph{cf.}~\cite{EGAII}*{corollaire~6.6.2}).
\blem\label{descent-ample-bundle}
Let $S$ be an integral Pr\"ufer scheme with fraction field $K$.
Let $L/K$ be a finite field extension and $T$ the normalization of $S$ in $L$.
For a smooth $S$-scheme $X$, denote its base change to $T$ by $X_T$.
\benumr 
\item\label{norm-map} There is norm map of invertible sheaves that induces a homomorphism of abelian groups
\[
 \Norm_f\colon \Pic(X_T)\ra \Pic(X) 
\]
such that for every invertible $\sO_X$-module $\sM$, we have $\Norm_f(f^{\ast}\sM)\simeq \sM^{\otimes d}$, where $d=[L:K]$.

\item\label{descend-ample} If $X_T$ has a $T$-ample invertible sheaf $\sL$, then $X$ has an $S$-ample invertible sheaf $\sE$.
\item\label{descend-quasi-affine} $X_T$ is $T$-quasi-affine if and only if $X$ is $S$-quasi-affine.
\eenum 
\elem 
\bpf 
%A standard limit argument reduces us to the case when $V$ has finite rank.
%Write $V$ as a filtered direct union of finite-rank valuation subdomains, taking integral closure commutes with this filtered direct limit.
%Then invertible sheaves on $X$ descend and by \cite{EGAIV3}*{lemme~8.10.5.2} so does the ampleness.
The assertion \ref{descend-quasi-affine} follows from \ref{descend-ample}.
We first prove the assertion \ref{norm-map} in two steps as follows. \\
\textbf{Step 1.} We claim that $f\colon X_T\ra X$ induces a multiplicative map of sheaves
\[
   \Norm_f\colon f_{\ast}\sO_{X_T}\ra \sO_{X}
\]
such that the composite $\sO_X\overset{f^{\sharp}}{\ra} \sO_{X_T}\overset{\Norm_f}{\lra} \sO_X$ sends $s$ to $s^d$, where $d=[L:K]$, and for any open $U\subset X$, if $r\in \GG(f^{-1}(U),\sO_{X_T})$ vanishes at $x\in f^{-1}(U)$, then $\Norm_f(r)=0$ at $f(x)$.

Since the claim is local and normalization commutes with smooth base change \SP{03GV}, we may assume that $X=\Spec A$ and $X_T=\Spec B$ for normal domains $A$ and $B$.
There is a norm map
\[
   \Norm_{F_B/F_A}\colon F_B^{\times}\ra F_A^{\times},\q l\mapsto \det(F_B\overset{l}{\ra} F_B),
\]
where $F_{A}$ and $F_B$ are fraction fields.
By \SP{0BIG} the characteristic polynomial $\Phi_l(T)$ of $l$ is a power of the minimal polynomial.
 If $l\in B$, then by \SP{00H7}, the constant term $\Norm_{F_B/F_A}(l)$ of $\Phi_l(T)$ is  in $A$. 
Consequently, we get a norm map $\Norm_{B/A}\ce \Norm_{F_B/F_A}|_B$ and its global version $\Norm_f$.
By construction, $\Norm_f\circ f^{\sharp}(s)=s^d$.
Let $b\in B$ be contained in the prime $\fp\subset A$. 
Then the constant term of the minimal polynomial is in $A\cap \fp$, so the norm of $b$ vanishes at $\fp$.
Hence, the claim holds.

\textbf{Step 2.}  By \Cref{one-dim-semilocal}, each fibre of the map $T\to S$ is a finite set of points. 
Therefore, $X_T\ra X$ is an integral morphism with finite fibres. 
Thus, for every invertible $\sO_{X_T}$-module $\sL$ and a point $x\in X$, by \SP{0F20}, there is $s\in \GG(X_T,\sL)$ that does not vanish at  $f^{-1}(x)$. 
In particular, there is an open neighborhood $U\subset X$ of $x$ such that $f^{-1}(U)\subset (X_T)_s$ so $\sL|_{f^{-1}(U)}$ is trivial.
Therefore, there is an open covering $X=\bigcup U_i$ such that every $\sL|_{f^{-1}(U_i)}$ is trivial.
Choose generating sections $s_i\in \GG(f^{-1}(U_i),\sL)$ and consider cocycles $\alpha_{ij}\in \sO^{\ast}(f^{-1}(U_i)\cap f^{-1}(U_j))$ determined by $s_i=\alpha_{ij}s_j$.
Then $\Norm_f(\alpha_{ij})$ form cocycles by the multiplicative property of $\Norm_f$ and defines an invertible $\sO_X$-module $\sE$.
By construction, the map $\Norm_f\colon \Pic(X_T)\ra \Pic(X)$ is additive map of Picard groups such that $\Norm_f(f^{\ast}\sM)=\sM^{\otimes d}$.

It remains to show \ref{descend-ample}.
By the first two steps, for the $T$-ample invertible sheaf $\sL$, its norm $\sE\ce\Norm_f(\sL)$ is an invertible $\sO_X$-module. 
Since $X$ is quasi-compact, it remains to show that $\sE$ is $S$-ample by checking:  every $x\in X$ has an affine open neighborhood $X_s$ for some $s\in \GG(X,\sE^{\otimes n})$.
For this, we choose an affine open $U\subset X$ containing $x$ and for $f^{-1}(x)\subset f^{-1}(U)$ exploit \cite{EGAII}*{corollaire~4.5.4} to obtain $t\in \GG(X_T,\sL^{\otimes k})$ such that $(X_T)_t$ is an affine open neighborhood of $f^{-1}(x)$.
Then the image $\tau \ce \Norm_f(t)$ is a section of $\sE^{\otimes k}$ such that $X_{\tau}$ is an open neighborhood of $x$, which is affine by \SP{05YU}, as desired.
\epf
\brem 
In \Cref{descent-ample-bundle}, if $L/K$ is a finite separable field extension, then $R$ is contained in a $V$-finite module, see \cite{Bou98}*{chapitre~V, \S1, \textnumero6, proposition~18}.
Roughly speaking, there is a basis $(w_1,\cdots,w_n)$ for $L/K$ such that $w_i\in R$, and the presence of a nondegenerate trace form $\mathrm{Tr}_{L/K}$ due to the separability of $L/K$ yields a dual base $(w_i^{\ast})_{i=1}^n$ such that $\sum_{i=1}^nVw_i\subset R\subset \sum_{i=1}^nVw_i^{\ast}\subset d^{-1}\sum_{i=1}^nVw_i$ (by leveraging the trace form), where $d$ is the discriminant of $L/K$ with respect to $(w_i)_{i=1}^n$.
However, the separable condition is superfluous when proving \Cref{descent-ample-bundle}.
In the Noetherian case, $R$ is a finite $V$-module.
This fails when $V$ is a nondiscretely valued, see Ostrowski's example in \cite{Rib99}*{6.3, Example~2}.
\erem

\blem \label{discon-degen-nonred-discon}
For a valuation ring $V$ of finite rank $n$ with spectrum $(S,\eta)$ and a connected scheme $X$ flat over $S$ whose every nongeneric fibre is reduced, every localization of $X$ at a nonclosed point of $S$ is connected.
In particular, the generic fibre $X_{\eta}$ is connected.
\elem
\bpf
Let $\fp\subset V$ be a prime of height $n-1$.
If $X_{\fp}$ is disconnected, then there is an idempotent $e\in \GG(X_{\fp}, \sO_{X_{\fp}})$.
It suffices to show that $e\in \GG(X,\sO_X)$.
For each affine open $U\ce\Spec R\subset X$ we have $e|_{U_{\fp}}\in R[\f{1}{\varpi}]$, where $\varpi\in \fm_V$ is an element of height $n$.
Hence, there is $r\in R$ and integer $k\geq 0$ such that $r/\varpi^k=e$ so we have $r^2=(e\varpi^k)^2=e\varpi^{2k}=\varpi^kr$.
If $k>0$, then $r$ is nilpotent in the closed fibre, contradicting with the reduced assumption. 
Therefore, $k=0$ so $e=r$ and we get $e\in \GG(X,\sO_{X})$.  
Replacing $X$ by $X_{\fp}$, the above argument reduces the rank $n$ to $n-1$, so the assertion follows.
\epf

\bthm\label{extend-extend} 
Let $S$ be an integral Pr\"ufer scheme with generic point $\eta$ and $G$ an $S$-smooth group.
For a $G$-homogeneous space $Y$, a $G^{\circ}$-stable open $X$ finitely presented over $S$, and an invertible sheaf $\sL_{\eta}$ on $X_{\eta}$,
\benumr
\item\label{basefree-power-extend-basefree} 
if $\sL_{\eta}$ is base-point-free, then it extends to an $S$-semiample invertible $\sO_X$-module $\sL$;
\item\label{genericample-extendample} if all nongeneric $S$-fibres of $X$ are connected, then any extension of an ample $\sL_{\eta}$ to $X$ is $S$-ample;
\item\label{separated-genericample-power-extend-ample} if  $S$ is a quasi-compact scheme whose closed subsets are finite sets\footnote{This happens when $S$ is normal, locally Noetherian of dimension one, or $S$ is a finite dimensional semilocal Pr\"ufer scheme.}, $X$ is $S$-separated and $\sL_{\eta}$ is ample, then there is a power $\sL_{\eta}^{\otimes n}$ that extends to an $S$-ample sheaf $\sM$ on $X$.
\eenum
\ethm
\bpf 
\ref{basefree-power-extend-basefree} 
Since $\sL_{\eta}$ is generated by its global sections, for each irreducible component $X_{\eta,i}$ of $X_{\eta}$, the subspace of $H^0(X_{\eta}, \sL_{\eta})$ that vanishes on $X_{\eta,i}$ is a proper subspace.
Because $k(\eta)$ is an infinite field, there is always a regular $k(\eta)$-linear combination $s_{\eta}$ of global sections of $\sL_{\eta}$.
Hence, the pair $(\sL_{\eta}, s_{\eta})$ determines an effective divisor $D_{\eta}$ such that $\sL_{\eta}\simeq \sO_{X_{\eta}}(D_{\eta})$ and $\Supp(D_{\eta})=\{s_{\eta}=0\}$.
By \Cref{Ramanujam--Samuel}, the schematic closure of $\{s_{\eta}=0\}$ in $X$ is a relative effective divisor $D$ and is $S$-flat.
Consequently, the invertible sheaf $\sL\ce \sO_{X}(D)$ is generated by $\mathbf{1}_{\sL}$ on the open $X\backslash \Supp(D)$, which contains all the maximal points of fibres of $X$.
It suffices to exploit \Cref{homo-ample-crit}\ref{semiample-crit} to conclude; note that after taking a power, $\sL$ is Picard-admissible with respect to the $G^{\circ}$-action by \Cref{C-criterion}.

\ref{genericample-extendample} For an extension $\sL$ of  $\sL_{\eta}$ to $X$, to show that $\sL$ is $S$-ample, by replacing $X$ with one of its irreducible components, we may assume that $X$ is integral and replace $S$ by the open image of $X$ such that $X\to S$ is surjective.
By spreading out and a limit argument \cite{EGAIV3}*{lemme~8.10.5.2}, we may assume that $S=\Spec V$ for a finite-rank valuation ring $V$. 
Since all nongeneric fibres of $X$ are connected, by \Cref{discon-degen-nonred-discon}, all $S$-fibres of $X$ are integral.
There is an $S$-affine open $U\subset X$ containing all the maximal points of $S$-fibres of $X$.
Hence, \Cref{affine-contain-ample-power-extend-ample} applies, yielding an invertible sheaf $\sM$ on $X$ extending a power $\sL_{\eta}^{\otimes n}$.
Since $\sL^{\otimes n}\otimes_{\sO_X}\sM^{-1}$ is trivial over $\eta$, by \Cref{gen-trivial}, it is the pullback of an invertible sheaf on $S$. 
Tensoring the ample sheaf $\sM$ to this pullback, we get the $S$-ampleness of  $\sL^{\otimes n}$ so of $\sL$.

\ref{separated-genericample-power-extend-ample}
There is a finite separable extension $k(\eta\pr)/k(\eta)$ provided by \cite{EGAIV2}*{corollaire~4.5.11} such that all the irreducible components of $X_{k(\eta\pr)}$ are geometrically irreducible. 
Let $T$ be the normalization of $S$ in $k(\eta\pr)$ with generic point $\eta\pr$.
If a power $\sL^{\otimes n}_{\eta\pr}$ extends to a $T$-ample invertible sheaf $\sM$ on $X_T$, then the norm map in  \Cref{descent-ample-bundle}  yields an $S$-ample invertible sheaf $\sE$ on $X$ with $\sE_{\eta}\simeq \sL_{\eta}^{n[L\colon K]}$ as a desired extension.
Hence, we may argue the assertion on $T$ and replace $X$ by an irreducible component so that we are reduced to the case when the generic fibre $X_{\eta}$ is geometrically irreducible. 

By spreading out \cite{EGAIV3}*{théorème~9.7.7}, there is an open subset $U\subset S$ such that every fibre of $X_U$ is irreducible.
Hence, we apply \ref{basefree-power-extend-basefree} and \ref{genericample-extendample} to $\sL_{\eta}$ to get an ample extension $\sL$ on $U$. 
We induct by showing that  $\sL$ extends over $U\cup \{z\}$ for an arbitrary maximal point $z\in Z$.
Let $U_z\subset S$ be an affine open neighborhood of $z$ such that $U_z\cap Z=\{z\}$. 
As $X$ is $S$-separated, by \Cref{finite-contain}, there is an affine open $W\subset X$ containing all the maximal points of $X_z$. 
Shrinking $U_z$ if necessary, we may assume that $W\surjects U_z$ is surjective.
For every $u\in U_z\backslash\{z\}$, the fibre $X_u$ is irreducible, hence the $U_z$-affine open subset $W$ contains all the maximal points $U_z$-fibres of $X$. 
By \Cref{affine-contain-ample-power-extend-ample}, there is a $U_z$-ample invertible sheaf $\sH$ on $X_{U_z}$ such that $\sH_{\eta}\simeq \sL^{\otimes k}_{\eta}$ for some $k>0$. 
The invertible sheaf $\sL_{U_z\cap U}^{\otimes k}\otimes \sH^{-1}_{U_z\cap U}$ is trivial on the generic fibre of $X_{U_z\cap U}$, hence by \Cref{general-ker-gen}, is the pullback of an invertible sheaf $\sN$ on $U_z\cap U$:
\[
\sL^{\otimes k}|_{X_{U_z\cap U}}\simeq \sH|_{U_z\cap U}\otimes \pi^{\ast}(\sN), 
\]
where $\pi\colon X_{U_z\cap U}\to U_z\cap U$ is the structural morphism.
Gluing $\sN$ with the trivial invertible sheaf on a small open neighborhood of $z$, we extend $\sN$ to an invertible $\sO_{U_z}$-module $\sN_{U_z}$. 
Consequently, $\sH_{X_{U_z\cap U}}\otimes \pi^{\ast}(\sN)$ extends to $\sH\otimes \pi^{\ast}(\sN_{U_z})$, which is an extension of $\sL^{\otimes k}_{X_{U_z\cap U}}$ to $X_{U_z}$.
It suffices to glue $\sL^{\otimes k}$ with $\sH\otimes \pi^{\ast}(\sN)$ to obtain the desired ample invertible sheaf on $U\cup U_z$. 
Since $Z\ce S\backslash U$ is a finite set, our induction stops in finite steps, so a power of $\sL_{\eta}$ extends to an ample invertible sheaf on $X$.
\qedhere

%Since $\sL_{\eta}$ is ample, there exists $n > 0$ and a regular section $s_\eta \in \Gamma(X_\eta, \sL_\eta^{\otimes n})$ such that the open subset $U_\eta \ce (X_\eta)_{s_\eta}$ is affine.
%Let $D$ be the schematic closure of the effective divisor $\mathrm{div}(s_\eta)$ in $X$.
%Since the property of being a relative effective divisor is local on the base, and $D$ restricted to any affine open subset of $S$ is the schematic closure of $\mathrm{div}(s_\eta)$ over a Pr\"ufer domain, \Cref{Ramanujam--Samuel} implies that $D$ is a relative effective divisor on $X$.
%We define the global extension $\sM \ce \sO_X(D)$. Note that $\sM_{\eta} \cong \sL_{\eta}^{\otimes n}$.
%
%We claim that $\sM$ is $S$-ample. 
%As ampleness is local on the base \SP{01VJ} and $S$ is tree-like partially ordered, by spreading out, we may assume that $S=\Spec V$ for a valuation ring $V$. 
%The open  $U \ce X \setminus \Supp(D)$ contains all the maximal points of fibres of $X$ (as $D$ is flat of pure codimension one).
%
%The pair $(X, \sM)$ satisfies the conditions of \ref{genericample-extendample}: $X_V$ is separated, $X_V \setminus \Supp(D_V)$ contains all maximal points of fibers, and the generic fiber is affine.
%By the same argument as in \ref{genericample-extendample} (specifically, using \Cref{descent-ample-bundle} to deduce ampleness via descent from a finite cover over a valuation base, or applying \Cref{cor-prufer-qproj}), we conclude that $\sM|_{X_V}$ is ample.
%Thus, $\sM$ is globally $S$-ample.
\epf
\bcor 
Let $S$ be an integral Pr\"ufer scheme with generic point $\eta$ and $G$ an $S$-smooth group.
For a $G$-homogeneous space $Y$, an open subscheme $X\subset Y$ with connected $S$-fibres, a closed subscheme $Z_{\eta}\subset X_{\eta}$ and its schematic closure $Z$ in $X$ with complementary open $U\ce X\backslash Z$, then
\[
   \text{$U$ is $S$-quasi-affine\q if and only if \q $U_{\eta}$ is $\eta$-quasi-affine.}
\]
\ecor
\bpf 
We only need to show the sufficiency.
Replacing $X$ by $G^{\circ}\cdot X$ and $Z_{\eta}$ by $(G^{\circ}\cdot X)\backslash U_{\eta}$, we may assume that $X$ is stable under $G^{\circ}$.
Let $D_{\eta}$ be an effective divisor such that the maximal points of $\Supp(D_{\eta})$ are the one-codimensional points of $Z_{\eta}$  in $X_{\eta}$.
Namely, there is an invertible $\sO_{X_{\eta}}$-module $\sL_{\eta}$ and a regular section $s_{D_{\eta}}\in \GG(X_{\eta},\sL_{\eta})$ such that $\Supp(D_{\eta})=V(s_D)\ce\{s_D=0\}$.
The schematic closure $\ov{V(s_D)}$ of $V(s_D)$ in $X$ by \Cref{Ramanujam--Samuel} determines an effective divisor $D$ on $X$. 
By \Cref{VI2.3}, the invertible $\sO_{X_{\eta}}$-module $\sO_{X_{\eta}}(D_{\eta})$ is ample, hence \Cref{extend-extend}\ref{genericample-extendample} implies that  $\sO_X(D)$ is $S$-ample.
Subsequently, $X\backslash \Supp(D)$ is $S$-quasi-affine, so  $\Supp(D)\subset Z$ implies that $U=X\backslash Z$ is $S$-quasi-affine.
\epf

\bcor 
For a connected scheme $S$, an $S$-smooth group scheme $G$ with connected $S$-fibres, a $G$-homogeneous space $X$ of finite type over $S$ with proper $S$-fibres, and an invertible $\sO_X$-module $\sL$, 
\[
   \text{if $\sL_s$ is ample for an $s\in S$,\q then $\sL$ is $S$-ample and $X$ is $S$-projective.}
\]

\ecor
\bpf 
For the subset $\cA \ce \{s\in S: \sL_s\;\text{is ample}\}$ of $S$, we prove that it is clopen such that $\sL|_{X_{\cA}}$ is ample and $X_{\cA}$ is $\cA$-projective.
For the openness, we show that each $s\in \cA$ has an open neighborhood $U(s)$ such that $X_{U(s)}$ is $U(s)$-projective and $\sL|_{X_{U(s)}}$ is $U(s)$-ample.
Note that $X\ra S$ is flat and locally of finite presentation due to \Cref{raygru}, by \cite{EGAIV2}*{corollaire~2.3.12}, it is universally submersive.
By \cite{SGA3Inew}*{exposé~VI$_{\text{B}}$, théorème~5.3}, $X$ is $S$-separated with proper geometrically connected fibres. 
Thus, a local criterion for properness \cite{EGAIV3}*{15.7.8} implies that $X$ is proper over $S$. 
Thanks to the $S$-properness of $X$, the ample aspect follows from \cite{EGAIV3}*{corollaire~9.6.4}.
It remains to show that $\cA$ is closed.
By the retrocompactness of $\cA\subset S$ and \cite{EGAII}*{proposition~7.2.2}, it suffices to show that $\cA$ is stable under specialization.
Let $y\in \cA$ nontrivially specialize to $y\pr\in S$ and let $x\in X$ be a point lying over $y$. 
For the morphism $X\ra S$, by \cite{EGAII}*{proposition~7.1.4}, there is a valuation ring $V$ with generic point $\eta$ and closed point $t$, a morphism $f\colon \Spec V\ra S$ such that $f(\eta)=y$ and $f(t)=y\pr$, and a rational $S$-map sending $\eta$ to $x$ such that $k(x)\simeq k(\eta)$.
Then, for the base change $X_V\ra \Spec V$ and the $\eta$-ample invertible $\sO_{X_{\eta}}$-module $\sL_{\eta}$, \Cref{extend-extend}\ref{genericample-extendample} implies that $\sL_{t}$ is ample.
Finally, the fpqc descent of ampleness \cite{EGAIV2}*{2.7.2}
 implies that $y\pr\in \cA$, so we conclude.
\epf

%\bibliographystyle{habbrv.bst}
%\bibliography{Cf1}

\end{document}